\documentclass[11pts]{article}

\usepackage{amssymb,latexsym, mathrsfs,amsmath}
\usepackage{graphicx}
\usepackage{verbatim}
\usepackage[all]{xy}
\usepackage{setspace}
\usepackage{amscd}

\newtheorem{theorem}{Theorem}[section]
\newtheorem{lemma}[theorem]{Lemma}
\newtheorem{proposition}[theorem]{Proposition}
\newtheorem{corollary}[theorem]{Corollary}
\newenvironment{definition}[1][Definition]
{\begin{trivlist} \item[\hskip \labelsep {\bfseries #1}]}
{\end{trivlist}}

\newenvironment{remark}[1][Remark]
{\begin{trivlist}
\item[\hskip \labelsep {\bfseries #1}]}
{\end{trivlist}}
\newenvironment{conjecture}[1][Conjecture]
{\begin{trivlist} \item[\hskip \labelsep {\bfseries #1}]}
{\end{trivlist}}
\newenvironment{proof}[1][Proof]
{\begin{trivlist}
\item[\hskip \labelsep {\bfseries #1}]}
{\end{trivlist}}
\newenvironment{acknowledgement}[1][Acknowledgement]
{\begin{trivlist}
\item[\hskip \labelsep {\bfseries #1}]}
{\end{trivlist}}
\newenvironment{affiliation}[1][]
{\begin{trivlist}
\item[\hskip \labelsep {\bfseries #1}]}
{\end{trivlist}}

\newcommand{\susp}{\Sigma\mkern-10.1mu/}

\DeclareMathSymbol{\N}{\mathbin}{AMSb}{"4E}
\DeclareMathSymbol{\Z}{\mathbin}{AMSb}{"5A}
\DeclareMathSymbol{\R}{\mathbin}{AMSb}{"52}
\DeclareMathSymbol{\Q}{\mathbin}{AMSb}{"51}
\DeclareMathSymbol{\I}{\mathbin}{AMSb}{"49}
\DeclareMathSymbol{\C}{\mathbin}{AMSb}{"43}

\def\P{{\mathbb P}}
\def\RP{{\mathbb R \mathbb P}}
\def\CP{{\mathbb C \mathbb P}}
\def\A{{\mathbb A}}

\def\S{{\mathbb S}}

\title{A homology and cohomology theory for real projective
varieties}
\author{Jyh-Haur Teh}

\begin{document}
\maketitle

\begin{abstract}
In this paper we develop homology and cohomology theories which play
the same role for real projective varieties that Lawson homology and
morphic cohomology play for projective varieties respectively. They
have nice properties such as the existence of long exact sequences,
the homotopy invariance, the Lawson suspension property, the
homotopy property for bundle projection, the splitting principle,
the cup product, the slant product and the natural transformations
to singular theories. The Friedlander-Lawson moving lemma is used to
prove a duality theorem between these two theories. This duality
theorem is compatible with the $\Z_2$-Poincar\'e duality for real
projective varieties with connected full real points.
\end{abstract}

\section{Introduction} The study of solving polynomial equations
dates back to the very beginning of mathematics. Finding general
algebraic solutions of a given equation was the original goal. This
goal was achieved for equations of degree 2, 3 and 4. But it was
proved by Abel and Galois that it was impossible for equations of
degree 5. Galois theory was created several decades later to study
some properties of the roots of equations. At the same time, people
started to consider the more complicated problem of solving
polynomial equations of more than one variable. The zero loci of
polynomial equations, which are called algebraic varieties, are
basic source of geometry and exemplify many important geometric
phenomena. In this paper, we study the properties of projective
varieties which are the zero loci of homogeneous polynomials in
projective spaces.

Algebraic cycles are some finite formal sum of irreducible
subvarieties with integral coefficients. The group $Z_p(X)$ of
$p$-cycles of a projective variety $X$ encodes many properties of
$X$. For many decades, since $Z_p(X)$ is a very large group in
general, quotients of $Z_p(X)$ were studied instead. For example,
the quotient of $Z_p(X)$ by rational equivalence is the Chow group
$CH_p(X)$ on which the intersection theory of algebraic varieties
can be built \cite{Fu}.

For a complex projective variety $X$, the group $Z_p(X)$ has
additional natural structure. According to the Chow theorem, there
is a canonical way, by means of Chow forms, to give $\mathscr{C}_{p,
d}(X)$, the set of effective $p$-cycles of degree $d$, the structure
of a projective variety. Thus $\mathscr{C}_{p, d}(X)$ has a
canonical analytic topology. Since $Z_p(X)$ is the group completion
of the monoid $\underset{d\geq 0}{\prod}\mathscr{C}_{p, d}(X)$, it
inherits a topology from this analytic structure which makes it a
topological group. From this point of view, methods from topology,
especially from homotopy theory, come into play. The title of the
Hirzebruch's book, ``Topological methods in algebraic geometry"
\cite{Hir}, perfectly describes the way we study projective
varieties via their cycle groups. The starting point of this
approach is the Lawson suspension theorem \cite{Lawson1}, which says
that there is a homotopy equivalence between $Z_p(X)$ and
$Z_{p+1}(\susp X)$, where $\susp X$ is the projective cone over $X$
(or more topologically the Thom space of the $\mathcal{O}(1)$ bundle
over $X$). It is then natural to define the Lawson homology group
$L_pH_n(X)=\pi_{n-2p}Z_p(X)$, the $(n-2p)$-th homotopy group of
$Z_p(X)$. When $p=0$, we have by the Dold-Thom theorem that
$L_0H_n(X)=H_n(X; \Z)$ (the singular homology of $X$). Thus we view
the Lawson homology as an enrichment of singular homology for
projective varieties. This point of view will be strengthened after
we develop a corresponding theory for real projective varieties and
extend some classical theorems, for example the Harnack-Thom theorem
\cite{Teh}, from singular homology.

In the past 15 years, Lawson, Friedlander, Mazur, Gabber,
Michelsohn, Lam, Lima-Filho, Walker and dos Santos have discovered
many properties of Lawson homology and have related it to other
theories like Chow groups, higher Chow groups, motivic cohomology
and algebraic $K$-theory. Notably, morphic cohomology was
established by Friedlander and Lawson, and a duality theorem between
morphic cohomology and Lawson homology was proved by using their
moving lemma \cite{FL1}, \cite{FL2}. It has been shown that Lawson
homology and morphic cohomology groups admit limit mixed Hodge
structure in \cite{FM}, \cite{FL1} and \cite{Walker}.

Algebraic topology also reaps the bounty of this harvest. Since
$$Z_0(A^n)=\frac{Z_0(\P^n)}{Z_0(\P^{n-1})}=K(\Z, 2n)$$ where
$K(\Z, 2n)$ is the Eilenberg-Mac Lane space, we are able to
represent many Eilenberg-Mac Lane spaces by more concrete algebraic
cycle spaces. For example, consider the natural embedding
$\mathscr{G}^q(\P^n)\subset \mathscr{C}^q(\P^n)$ of the Grassmannian
of codimension-q planes into the limit space
$\mathscr{C}^q(\P^n)=\lim_{d\to \infty}\mathscr{C}^q_d(\P^n)$ of
codimension-q cycles of $\P^n$. Letting $n$ go to infinity, we get a
map
$$BU_q \overset{c}{\longrightarrow} K(\Z, 2)\times K(\Z, 4) \times
\cdots \times K(\Z, 2q)$$ where $BU_q=\lim_{n\to
\infty}\mathscr{G}^q(\P^n)$ is the classifying space for the unitary
group $U_q$. A beautiful theorem by Lawson and Michelsohn \cite{LM}
says that for each $q\geq 1$, this map induces an isomorphism
$$\Z\cong \pi_{2q}(BU) \overset{c_*}{\longrightarrow}
\pi_{2q}(K(\Z, ev))\cong \Z$$ which is multiplication by $(q-1)!$
where $K(\Z, ev)=\prod^{\infty}_{i=1}K(\Z, 2i)$. Our next example is
from morphic cohomology. \begin{itemize}
    \item There is a join pairing in algebraic cocycle groups which
    induces the cup product in singular cohomology.
    \item The Chern classes in morphic cohomology map to the Chern
    classes in singular cohomology.
    \item The cocycle group $Z^m(X)$ of a $m$-dimensional smooth
    projective manifold is homotopy equivalent to the mapping space
    $Map(X, K(\Z, 2m))$.
\end{itemize}

The Fundamental Theorem of Algebra says that a degree $d$ polynomial
in one variable over $\C$ has $d$ zeros in complex plane counting
multiplicities. A polynomial $f$ of degree $d$ over $\R$ needs not
have $d$ real zeros, but if we take $Z(f)$ to be the set of $d$
zeros of $f$, and $Z(f)^{av}$ to be the set of nonreal zeros of $f$
counting multiplicities, we find that if we define
$R(f)=Z(f)-Z(f)^{av}$, since nonreal zeros appear in conjugate
pairs, the cardinality $|R(f)|$ of $R(f)$ and the degree $d$ have
the following relation:
$$|R(f)|\equiv d \mbox{ mod } 2$$
This is the Reduced Real Fundamental Theorem of Algebra.

The principle is that everything in the complex world has a
counterpart in the real world with $\Z_2$-coefficients. Over $\C$,
the utility of (co)homology theory with $\Z$-coefficients is related
to the fact that the zeros of any polynomial can be counted. But
over $\R$, we can only count the zeros modulo 2. The bridge from the
complex world to the real world passes to a quotient of the set of
elements invariant under conjugation. The reduced real cycle group
is defined as
$$R_p(X)=\frac{Z_p(X)_{\R}}{Z_p(X)^{av}}$$
where $X$ is a real projective variety, $Z_p(X)_{\R}$ is the set of
$p$-cycles which are invariant under conjugation and $Z_p(X)^{av}$
are cycles of the form $c+\bar{c}$. This group first appeared in the
thesis of Lam where he proved the Lawson suspension theorem for
reduced real cycle groups. We develop this idea further to get a
homology-like theory called reduced real Lawson homology and a
cohomology-like theory called reduced real morphic cohomology. In
section 2, we define the topology on some cycle groups. In section
3, we develop relative theories of reduced real Lawson homology and
reduced real morphic cohomology, and get some long exact sequences.
The functoriality properties of reduced real theories are considered
in section 4. We show in section 5 that the reduced real morphic
cohomology has many nice properties: the homotopy invariance, the
splitting principle, the suspension property, and the homotopy
property. In section 6, by using the moving lemma of Friedlander and
Lawson, we prove a duality theorem between reduced real Lawson
homology and reduced real morphic cohomology. As in morphic
cohomology, the join pairing of cycles induces a cup product in the
reduced real morphic cohomology and through a natural tranformation
to singular cohomology with $\Z_2$-coefficients it becomes the usual
cup product. The natural transformations and the $s$-maps in these
two theories define a filtration in the singular homology and
cohomology with $\Z_2$-coefficients respectively. Under the natural
transformations constructed in section 7, we show that our duality
is compatible with the usual Poincar\'e duality with
$\Z_2$-coefficients. In appendix, we show that all cycle groups that
we deal with in this paper are CW-complexes.

\section{Topology Of Reduced Real Cycle And Cocycle Spaces}
Throughout this paper, $X, Y$ are projective varieties of dimension
$m$ and $n$ respectively and $X$ is normal. An $r$-cycle on $X$ is a
linear combination of irreducible subvarieties of the same dimension
$r$ in $X$ with $\Z$-coefficients. Let $\mathscr{C}_{r, d}(X)$ be
the set of degree $d$ $r$-cycles of $X$ with positive coefficients.
By Chow theorem, $\mathscr{C}_{r, d}(X)$ has the structure of a
projective variety. Thus we may consider $\mathscr{C}_{r, d}(X)$ as
a complex projective variety with the analytic topology. Let $K_{r,
d}(X)=\underset{d_1+d_2\leq d}{\prod}(\mathscr{C}_{r, d_1}(X)\times
\mathscr{C}_{r, d_2}(X))/\sim$ where the equivalence relation $\sim$
is given by $(a_1, b_1)\sim (a_2, b_2)$ if and only if
$a_1+b_2=a_2+b_1$. $K_{r, d}(X)$ inherits topology from this
quotient which makes $K_{r, d}(X)$ a compact Hausdorff space. And
from the filtration,
$$K_{r, 0}(X) \subset K_{r, 1}(X) \subset \cdots =Z_r(X),$$
we give the space of $r$-cycles $Z_r(X)$ on $X$ the  weak topology
which means that a subset $C \subset Z_r(X)$ is closed if and only
if $C\cap K_{r, d}(X)$ is closed for all $d$. This topology makes
$Z_r(X)$ a topological abelian group. In general, let
$$A_1\subset A_2 \subset A_3 \subset \cdots  =A$$ be a chain of
closed inclusions of topological spaces. We define the topology on
$A$ by declaring a subset $C\subset A$ to be closed if and only if
its intersection $C\cap A_i$ is closed for all $i\geq 1$. This
topology is called the weak topology with respect to the subspaces.

\begin{definition}
Let $\mathscr{C}_r(Y)(X)$ be the topological submonoid of
$\mathscr{C}_{r+m}(X\times Y)$ consisting of elements which are
equidimensional of relative dimension $r$ over $X$. Let $Z_r(Y)(X)$
be the naive group completion $[\mathscr{C}_r(Y)(X)]^+$ of
$\mathscr{C}_r(Y)(X)$ with the quotient topology. We call
$Z_r(Y)(X)$ the dimension $r$ cocycle group of $X$ with values in
$Y$. Let $\mathscr{M}or(X, \mathscr{C}_r(Y))=\coprod_{d\geq
0}^{\infty}\mathscr{M}or(X, \mathscr{C}_{r, d}(Y))$ be the space of
all morphisms from $X$ to $\mathscr{C}_r(Y)$ with the compact-open
topology. Let $\mathscr{M}or(X, Z_r(Y))=[\mathscr{M}or(X,
\mathscr{C}_r(Y))]^+$ be the naive group completion of
$\mathscr{M}or(X, \mathscr{C}_r(Y))$.
\end{definition}

By the graphing construction of Friedlander and Lawson \cite[Theorem
1.2]{FL3}, there is a topological group isomorphism
$$\mathcal{G}: \mathscr{M}or(X, Z_r(Y))\rightarrow
Z_r(Y)(X)$$

\begin{definition}
A real projective variety $V\subset \CP^n$ is a complex projective
variety which is invariant under the conjugation of $\CP^n$. A
subvariety $V\subset X\times Y$ is real if and only if
$\overline{V}=V$ where $\overline{V}=\{(\overline{x},
\overline{y})|(x, y)\in V\}$. A cycle on $V$ is real if all its
components are real subvarieties of $V$. A cycle $c$ is averaged if
$c=a+\overline{a}$ for some cycle $a$.
\end{definition}

We will write the conjugate of a point $x\in X$ as $\overline{x}$.
The conjugation induces a map on cycle groups. For $f\in Z_r(Y)(X)$,
define $\overline{f}(x)=\overline{f(\overline{x})}$. This notation
is convenient when we think of $f$ as a map from $X$ to the cycle
space of $Y$.

\begin{proposition}
For $f\in Z_r(Y)(X)$, $f$ is real if and only if $\overline{f}=f$.
\end{proposition}

\begin{proof}
Since $f$ is real, the set $\{(x, f(x))|x\in X\}$ is equal to the
set $\{(\overline{x}, \overline{f(x)})|x\in X\}$. The point
$(\overline{x}, f(\overline{x}))$ is in this set and $f$ is a
function, hence $f(\overline{x})=\overline{f(x)}$. Therefore
$\overline{f}=f$. Another direction is trivial.
\end{proof}

\begin{definition}
The filtration $$K_{r, 0}(X) \subset K_{r, 1}(X) \subset K_{r, 2}(X)
\subset \cdots =Z_r(X)$$ is called the canonical $r$-filtration of
$X$ and if $X$ is a real projective variety,
$$K_{r, 0}(X)_{\R} \subset K_{r, 1}(X)_{\R} \subset K_{r, 2}(X)_{\R} \subset \cdots
=Z_r(X)_{\R}$$ is called the canonical real $r$-filtration of $X$
where $K_{r, i}(X)_{\R}$ is the subset of real cycles in $K_{r,
i}(X)$ and
$$K_{r, 0}(X)^{av} \subset K_{r, 1}(X) \subset K_{r, 2}(X)^{av} \subset \cdots
=Z_r(X)^{av}$$ is called the canonical averaged $r$-filtration of
$X$ where $K_{r, i}(X)^{av}$ is the subset of averaged cycles in
$K_{r, i}(X)$. If a filtration is defined by a sequence of compact
sets, this filtration is called a compactly filtered filtration.
\end{definition}

\begin{definition}
Let $Z_r(X)_{\R}$ be the subgroup of $Z_r(X)$ consisting of real
cycles and let $Z_r(X)^{av}=\{V+\overline{V}|V\in Z_r(X)\}$ be the
averaged cycle group. We give $Z_r(X)_{\R}$ and $Z_r(X)^{av}$ the
subspace topology of $Z_r(X)$ which are same as the weak topology
defined by the canonical real $r$-filtration and the canonical
averaged $r$-filtration respectively.
\end{definition}

\begin{proposition}
Suppose that $Y$ is a subvariety of a projective variety $X$. Then
$Z_r(Y)$ is a closed subgroup of $Z_r(X)$. If $Y\subset X$ are real
projective varieties, then $Z_r(Y)_{\R}, Z_r(Y)^{av}$ are closed
subgroups of $Z_r(X)_{\R}$ and $Z_r(X)^{av}$ respectively.
\end{proposition}

\begin{proof}
This is because $Z_r(Y)\cap K_{r, d}(X)=K_{r, d}(Y)$ which is closed
in $K_{r, d}(X).$ Similarly, the closedness follows from the fact
$Z_r(Y)_{\R}\cap K_{r, d}(X)_{\R}=K_{r, d}(Y)_{\R}$ and
$Z_r(Y)^{av}\cap K_{r, d}(X)^{av}=K_{r, d}(Y)^{av}$.
\end{proof}

Let $\mathscr{C}_{r, d}(Y)(X)$ be the subspace of $\mathscr{C}_{r+m,
d}(X\times Y)$ consisting of cycles which are equidimensional over
$X$. It is shown in \cite[Lemma 1.4]{FL1} that $\mathscr{C}_{r,
d}(Y)(X)$ is a Zariski open set of $\mathscr{C}_{r+m, d}(X\times
Y)$. Let $K_{r, d}(Y)(X)=\coprod_{d_1+d_2\leq d}\mathscr{C}_{r,
d_1}(Y)(X)\times \mathscr{C}_{r, d_2}(Y)(X)/\sim$ with the quotient
topology, and $K_{r, d}(Y)(X)_{\R}, K_{r, d}(Y)(X)^{av}$ the
subspaces of $K_{r, d}(Y)(X)$ consisting of elements which can be
represented by real cycles and averaged cycles respectively. Then we
may formed the corresponding canonical filtration, canonical real
filtration and canonical averaged filtration.

\begin{definition}
Suppose that $X, Y$ are real projective varieties. Let
$$Z_r(Y)(X)_{\R}=Z_r(Y)(X)\cap Z_{r+m}(X\times Y)_{\R}$$
$$Z_r(Y)(X)^{av}=Z_r(Y)(X)\cap Z_{r+m}(X\times Y)^{av}$$
with the subspace topology of $Z_r(Y)(X)$. Since the topology of
$Z_r(Y)(X)$ is given by the weak topology from the canonical
filtration, the topology of $Z_r(Y)(X)_{\R}$ and $Z_r(Y)(X)^{av}$
are given by the weak topology of the canonical real and averaged
filtration respectively. \label{real weak topology}
\end{definition}

\begin{definition}
A filtration of a topological space $T$ by a sequence of subspaces
$$T_0\subset T_1\subset \cdots \subset T_j \cdots =T$$
is said to be locally compact if for any compact subset $K\subset
T$, there exists some $e\geq 0$ such that $K\subset T_e$.
\end{definition}

\begin{lemma}
Suppose that $X$ is a Hausdorff space and the topology on $X$ is
given by the filtration
$$T_0\subset T_1\subset \cdots \subset T_j \cdots=X$$ by the weak
topology. Then this filtration is a locally compact filtration.
\label{compact subset}
\end{lemma}

\begin{proof}
Let $K\subset X$ be a compact subset. Assume that $K$ is not
contained in any $T_e$, then we can find a sequence $\{x_i\}_{i\in
I}$ where $x_i\in (T_i-T_{i-1})\cap K$ and $I$ is a sequence of
distinct integers which goes to infinity. Let $U_k=X-\{x_i\}_{i\in
I, i\neq k}$. Then $U_k\cap T_j=T_j-\{x_i|i\leq j, i\neq k\}$ which
is open for all $j$. Thus $U_k$ is an open subset of $X$ and
$K\subset \cup_{i\in I}U_k$ but we are unable to find a finite
subcover which contradicts to the compactness of $K$. So $K$ has to
be contained in some $T_e$.
\end{proof}

From this lemma, we know that all of our canonical filtrations are
locally compact.

\begin{proposition}
For real projective varieties $X, Y$, $Z_r(Y)(X)_{\R}$ and
$Z_r(Y)(X)^{av}$ are closed subgroups of $Z_r(Y)(X)$. In particular,
$Z_r(Y)_{\R}$ and $Z_r(Y)^{av}$ are closed subgroups of $Z_r(Y)$.
\label{closedness}
\end{proposition}

\begin{proof}
Define $\psi: Z_r(Y)(X) \longrightarrow Z_r(Y)(X)$ by
$\psi(f)=\overline{f}-f$. Since $\psi$ is continuous and
$Z_r(Y)(X)_{\R}$ is the kernel of $\psi$, $Z_r(Y)(X)_{\R}$ is
closed. Suppose that $\{f_n+\overline{f_n}\}$ is a sequence in
$Z_r(Y)(X)^{av}$ which converges to $c$. Since the canonical
averaged $r$-filtration of $X$ is locally compact, by Lemma
\ref{compact subset}, $A=\{f_n+\overline{f_n}\}\cup \{c\}$ is
contained in $K_{r, d}(Y)(X)^{av}$ for some $d>0$. Therefore there
exist $g_n\in K_{r, d}(Y)(X)^{av}$ such that
$g_n+\overline{g_n}=f_n+\overline{f_n}$ for all $n$. Since the
sequence $\{g_n\}\subset K_{r+m, d}(X\times Y)$ and $K_{r+m,
d}(X\times Y)$ is compact, hence $\{g_n\}$ has a convergent
subsequence, which we denote by $\{g_{n_i}\}$. Suppose that
$\{g_{n_i}\}$ converges to a cycle $g$. Then
$\{\overline{g_{n_i}}\}$ converges to $\overline{g}$. Therefore
$g_{n_i}+\overline{g_{n_i}}$ converges to $g+\overline{g}$ which
implies $c=g+\overline{g}$. Thus $Z_r(Y)(X)^{av}$ is closed. Take
$X$ to be a point, it follows immediately that $Z_r(Y)_{\R}$ and
$Z_r(Y)^{av}$ are closed in $Z_r(Y)$.
\end{proof}

\begin{definition}
Suppose that $X, Y$ are real projective varieties. We define the
reduced real cycle group to be
$$R_r(X)=\frac{Z_r(X)_{\R}}{Z_r(X)^{av}},$$
and the reduced real $Y$-valued cocycle group to be
$$R_r(Y)(X)=\frac{Z_r(Y)(X)_{\R}}{Z_r(Y)(X)^{av}}$$
Furthermore, if $Y$ is a real subvariety of $X$. We define $$Z_r(X,
Y)^{av}=\frac{Z_r(X)^{av}}{Z_r(Y)^{av}},$$ and
$$Z_r(X, Y)_{\R}=\frac{Z_r(X)_{\R}}{Z_r(Y)_{\R}}.$$ All these groups
are enrolled with quotient topology.
\end{definition}

A topological group $G$ is a group which is also a Hausdorff space
and for $(g, h)\in G\times G$, the product $(g, h)\longrightarrow
gh^{-1}$ is a continuous map. From the Proposition above, we have
the following result.

\begin{corollary}
For real projective varieties $X, Y$, $R_r(Y)(X)$ and $R_r(X)$ are
topological abelian groups.
\end{corollary}

\begin{proof}
This is from the general fact that the quotient of a topological
group by any of its closed normal subgroups is a topological group.
\end{proof}

The homotopy types of $Z_r(\P^n)_{\R}$ and $Z_r(\P^n)^{av}$ were
computed in \cite{LLM1} which are quite complicated but the homotopy
type of $R_r(\P^n)$ is much simpler. In the following we will see
that the reduced real cycle groups are closely related to the
singular homology with $\Z_2$-coefficients of the real points. These
are some of the reasons that we work on the reduced real cycle
groups.

There are two natural ways to put topology on $R_r(X)$. One is the
weak topology from the filtration induced by the canonical real
filtration on $R_r(X)$ another one is the quotient topology. We show
that these topologies coincide on $R_r(X)$.

\begin{proposition}
Suppose that $G$ is a topological group and
$$K_1 \subset K_2 \subset \cdots
=G$$ is a filtration filtered by compact subsets of $G$ which
generates the topology of $G$. Let $H$ be a closed normal subgroup
of $G$ and $q: G\longrightarrow G/H$ be the quotient map. Denote the
restriction of $q$ to $K_k$ by $q_k$ and $M_k=q(K_k)$. We define a
topology on $M_k$ by making $q_k$ a quotient map, for each $k$. Then
\begin{enumerate}
\item $M_k$ is a subspace of $M_{k+1}$ for all $k$.

\item the weak topology of $G/H$ defined by the filtration
$$M_1 \subset M_2 \subset \cdots =G/H$$
coincides with the quotient topology of $G/H$.
\end{enumerate}
\end{proposition}

\begin{proof}
\begin{enumerate}
\item If $C\subset M_{k+1}$ is a closed subset, since $M_{k+1}$ is
compact, $C$ is also compact. From the commutative diagram,
$$\xymatrix{ K_k \ar[d]^{q_k} \ar[r]^{i} & K_{k+1} \ar[d]^{q_{k+1}} & \\
M_k \ar[r]^{\bar{i}} & M_{k+1}}$$ we have $\bar{i}^{-1}(C)=q_k
i^{-1} q^{-1}_{k+1}(C)$ which is compact thus closed. So $\bar{i}$
is continuous. Since $M_k$ is compact and $\bar{i}$ is injective, it
is an embedding.

\item Let $C$ be a closed subset of $G/H$ under the quotient
topology. Then $q^{-1}(C)$ is closed which means $q^{-1}(C)\cap K_k$
is closed for all $k$. $q_k^{-1}(C \cap M_k)=q^{-1}(C)\cap K_k $ is
closed, so $C\cap M_k$ is closed for all $k$ thus $C$ is closed
under the weak topology. On the other hand, assume that $C$ is a
closed subset of $G/H$ under the weak topology, that is, $C\cap M_k$
is closed for all $k$. $q^{-1}(C)\cap K_k=q^{-1}(C \cap M_k)\cap
K_k=q_k^{-1}(C\cap M_k)$ which is closed for all $k$, thus $C$ is a
closed subset in the quotient topology.
\end{enumerate}
\end{proof}

\begin{corollary}
Suppose that $X$ is a real projective variety and $Y$ is a real
subvariety of $X$. Consider the canonical real and averaged
filtrations of $X$:
$$K_{r, 1}(X)_{\R} \subset K_{r, 2}(X)_{\R} \subset \cdots
=Z_r(X)_{\R}$$
$$K_{r, 1}(X)^{av} \subset K_{r, 2}(X)^{av} \subset \cdots
=Z_r(X)^{av}$$ and the quotient maps $q_1:
Z_r(X)_{\R}\longrightarrow R_r(X)$, $q_2: Z_r(X)_{\R}
\longrightarrow Z_r(X, Y)_{\R}$, $q_3: Z_r(X)^{av} \longrightarrow
Z_r(X, Y)^{av}$. Then
\begin{enumerate}
\item for all $k$, $q_1(K_{r, k}(X)_{\R})$, $q_2(K_{r,
k}(X)_{\R})$ and $q_3(K_{r, k}(X)^{av})$ are subspaces of $q_1(K_{r,
k+1}(X)_{\R})$, $q_2(K_{r, k+1}(X)_{\R})$ and $q_3(K_{r,
k+1}(X)^{av})$ respectively.

\item the weak topology of $R_r(X), Z_r(X, Y)_{\R}, Z_r(X,
Y)^{av}$ induced from the filtrations above coincide with the
quotient topology on them.
\end{enumerate}
\label{weak topology coincide}
\end{corollary}

\begin{proposition}
For a real projective variety $X$, let $ReX$ denote the set of real
points of $X$ and let $Z_0(ReX)$ denote the subgroup of $Z_0(X)$
generated by real points of $X$. Then
\begin{enumerate}
\item $Z_0(ReX)$ is a closed subgroup of $Z_0(X)_{\R}$.

\item $R_0(X)$ is isomorphic as a topological group to
$\frac{Z_0(ReX)}{2Z_0(ReX)}$.
\end{enumerate}
\label{real points isomorphic}
\end{proposition}

\begin{proof}
\begin{enumerate}
\item Since $X$ is compact, the quotient map from $X\times \cdots
\times X$(k-times) to the symmetric product $SP^k(X)$ is a closed
map. $ReX\times \cdots \times ReX$ is a closed subset of $X\times
\cdots \times X$ thus $SP^k(ReX)$ is a closed subset of $SP^k(X)$.
Consider the following filtrations:
$$A_0\subset A_1\subset A_2 \subset \cdots =Z_0(ReX)$$ where
$A_n=\coprod^n_{k=0}SP^k(ReX)\times SP^k(ReX)/ \sim$ and
$$B_0\subset B_1\subset B_2 \subset \cdots =Z_0(X)$$ where
$B_n=\coprod^n_{k=0}SP^k(X)\times SP^k(X)/ \sim$ and $\sim$ is the
equivalence relation coming from the naive group completion. The
topology of $Z_0(ReX)$ and $Z_0(X)$ are defined by these two
filtrations respectively. Since $SP^k(ReX)$ and $SP^k(ReX)$ are
compact, $A_n$ is closed in $B_n$. Observe that $Z_0(ReX)\cap
B_n=A_n$, thus $Z_0(ReX)$ is a closed subgroup of $Z_0(X)$ and hence
a closed subgroup of $Z_0(X)_{\R}$.

\item Let $q: Z_0(X)_{\R} \longrightarrow R_0(X)$ be the quotient
map and $i: Z_0(ReX)\hookrightarrow Z_0(X)_{\R}$ be the inclusion
map. Then $q\circ i: Z_0(ReX) \longrightarrow R_0(X)$ is a
continuous map and $2Z_0(ReX)$ is in the kernel, thus we get a
continuous map $\psi: \frac{Z_0(ReX)}{2Z_0(ReX)} \longrightarrow
R_0(X)$. Since each class in $R_0(X)$ can be represented uniquely by
cycle of real points modulo 2, $\psi$ is bijective. Let
$$K_0\subset K_1 \subset K_2 \subset \cdots =Z_0(X)_{\R}$$ be the canonical
real filtration and $T_i=K_i\cap Z_0(ReX)$. $Z_0(ReX)$ is a closed
subset of $Z_0(X)_{\R}$ and $K_i$ is compact, thus $T_i$ is compact
for all $i$. Since $\psi$ is bijective,
$\psi(T_i+2Z_0(ReX))=K_i+Z_0(X)^{av}$ for each $i\geq 0$. For a
closed $C\subset \frac{Z_0(ReX)}{2Z_0(ReX)}$, $\psi(C)\cap
(K_i+Z_0(X)^{av})=\psi(C\cap T_i+2Z_0(ReX))$ which is closed, thus
$\psi^{-1}$ is continuous. Then it is easy to see that $\psi$ is a
topological group isomorphism.
\end{enumerate}
\end{proof}

\begin{definition}
Let $X$ be a real projective variety. For any $f\in Z_p(X)$, let
$f=\sum_{i\in I}n_iV_i$ be in the reduced form, i.e., $V_i\neq V_j$
if $i\neq j$. Let
$$RP(f)=\sum_{i \in I, \overline{V_i}=V_i}n_i V_i$$ which is called
the real part of $f$. Let $$J=\{i \in I|V_i \mbox{ is not real and }
\overline{V_i} \mbox{ is also a component of } f \}$$ and for $i\in
J$, let $m_i$ be the maximum value of the coefficients of $V_i$ and
$\overline{V_i}$. Define the  averaged part to be
$$AP(f)=\sum_{i\in J}m_i(V_i+\overline{V}_i)$$ and the
imaginary part to be
$$IP(f)=f-RP(f)-AP(f).$$
\label{real part}
\end{definition}

Then $f\in Z_p(X)_{\R}$ if and only if $IP(f)=0$.

\begin{proposition}
Suppose that $X, Y, Z$ are real projective varieties and $Y$ is a
subvariety of $Z$, then
\begin{enumerate}
\item the inclusion map $i:Z_r(Y)(X) \longrightarrow Z_r(Z)(X)$ is a
closed embedding.

\item the inclusion map $i: Z_r(Y)(X) \longrightarrow Z_r(Z)(X)$
induces a closed embedding $\bar{i}: R_r(Y)(X) \hookrightarrow
R_r(Z)(X)$ for $r\geq 0$.

\item the inclusion map $i: Z_r(Z)^{av}\longrightarrow
Z_r(Z)_{\R}$ induces a closed embedding $\bar{i}: Z_r(Z, Y)^{av}
\hookrightarrow Z_r(Z, Y)_{\R}$ for $r\geq 0$.
\end{enumerate}
\label{subspace}
\end{proposition}

\begin{proof}
\begin{enumerate}
\item Let
$$G_{r, d}(Y)(X)=\coprod_{d_1+d_2\leq d}\mathscr{M}or(X,
\mathscr{C}_{r, d_1}(Y))\times \mathscr{M}or(X, \mathscr{C}_{r,
d_2}(Y)),$$ $$G_{r, d}(Z)(X)=\coprod_{d_1+d_2\leq d}\mathscr{M}or(X,
\mathscr{C}_{r, d_1}(Z))\times \mathscr{M}or(X, \mathscr{C}_{r,
d_2}(Z)).$$ Let $p:\mathscr{M}or(X, \mathscr{C}_r(Y))\times
\mathscr{M}or(X, \mathscr{C}_r(Y)) \rightarrow Z_r(Y)(X)$ and
$q:\mathscr{M}or(X, \mathscr{C}_r(Z))\times \mathscr{M}or(X,
\mathscr{C}_r(Z)) \rightarrow Z_r(Z)(X)$ be the quotient maps. Let
$G'_{r, d}(Y)(X)=p(G_{r, d}(Y)(X)), G'_{r, d}(Z)(X)=q(G_{r,
d}(Z)(X))$ enrolled with the quotient topology. Then by the graphing
construction of Friedlander and Lawson, the topology of $Z_r(Y)(X)$
and $Z_r(Z)(X)$ are same as the weak topologies induced by the
filtrations formed by $G'_{r,d}(Y)(X)$ and $G'_{r, d}(Z)(X)$
respectively. The topology of $\mathscr{M}or(X, \mathscr{C}_{r,
d}(Z))$ is the compact-open topology. Since $\mathscr{C}_{r, d}(Y)$
is a closed subspace of $\mathscr{C}_{r, d}(Z)$, it is easy to see
that $\mathscr{M}or(X, \mathscr{C}_{r, d}(Y))$ is a closed subspace
of $\mathscr{M}or(X, \mathscr{C}_{r, d}(Z))$. Consider $G'_{r,
d}(Y)(X)$ as a subset of $G'_{r, d}(Z)(X)$, then $p^{-1}(G'_{r,
d}(Y)(X))=G_{r, d}(Y)(X)$ which is a closed subset of $G_{r,
d}(Z)(X)$. Hence $G'_{r, d}(Y)(X)$ is a closed subset of $G'_{r,
d}(Z)(X)$. Therefore $Z_r(Y)(X)$ is a closed subgroup of
$Z_r(Z)(X)$.

\item Since the inclusion $Z_r(Y)(X) \overset{i}{\hookrightarrow}
Z_r(Z)(X)$ is an embedding, the restriction
$$i: Z_r(Y)(X)_{\R} \hookrightarrow
Z_r(Z)(X)_{\R}$$ is also an embedding. Since $Z_r(Y)(X)^{av} \subset
Z_r(Z)(X)^{av}$, $i$ induces a map
$$\bar{i}: R_r(Y)(X) \longrightarrow R_r(Z)(X).$$
If $f\in Z_r(Z)(X)^{av}\cap Z_r(Y)(X)_{\R}$, then $RP(f)=2g$ for
some $g\in Z_r(Y)(X)_{\R}$. Thus $f=2g+AP(f) \in Z_r(Y)(X)^{av}$ and
therefore $\bar{i}$ is injective.

Let $A_0 \subset A_1 \subset A_2 \subset \cdots =Z_r(Y)(X)_{\R}$,
$B_0 \subset B_1 \subset B_2 \subset \cdots =Z_r(Z)(X)_{\R}$ be the
canonical real filtrations. Let $q_1, q_2$ be the quotient maps from
$Z_r(Y)(X)_{\R}$, $Z_r(Z)(X)_{\R}$ to $R_r(Y)(X)$ and $R_r(Z)(X)$
respectively. Let $C_k=q_1(A_k), D_k=q_2(B_k)$ for all $k$. Since
$i(A_k) \subset B_k$, $\bar{i}(C_k) \subset D_k$ and from the
definition of canonical filtrations, $\bar{i}^{-1}(D_k)\subset C_k$.
For any closed subset $W$ of $R_r(Y)(X)$, $W\cap C_k$ is compact and
by the injectivity of $\bar{i}$, $\bar{i}(W\cap C_k)=\bar{i}(W)\cap
\bar{i}(C_k)=\bar{i}(W)\cap D_k$ which is compact and thus closed.

\item Since the inclusion $i: Z_r(Z)^{av} \hookrightarrow
Z_r(Z)_{\R}$ is an embedding and $i(Z_r(Y)^{av}) \subset
Z_r(Y)_{\R}$ , it induces a map $\bar{i}: Z_r(Z, Y)^{av}
\hookrightarrow Z_r(Z, Y)_{\R}$. If
$\bar{i}(f+Z_r(Y)^{av})=f+Z_r(Y)_{\R}=Z_r(Y)_{\R}$, then $f\in
Z_r(Y)_{\R}$. This implies $RP(f)\in Z_r(Y)_{\R}$. But $f$ is also
an averaged cycle, so $RP(f)=2g$ for some $g\in Z_r(Y)_{\R}$ and
$f=2g+AP(f)\in Z_r(Y)^{av}$. Hence, $\bar{i}$ is injective. For the
rest, it is similar to the argument above.
\end{enumerate}
\end{proof}

From this Proposition, when we need to take quotient, we will abuse
of notation and write $\frac{R_r(Z)(X)}{R_r(Y)(X)}$ and
$\frac{Z_r(Z, Y)_{\R}}{Z_r(Z, Y)^{av}}$ for the quotient of
$R_r(Z)(X), Z_r(Z, Y)_{\R}$ by the images of $\bar{i}$ in
$R_r(Y)(X)$ and $Z_r(Z, Y)^{av}$ respectively.

\section{Relative Theory And Long Exact Sequences}
In \cite[Proposition 4.1]{Teh}, we develop a technique to produce a
long exact sequence of homotopy groups from the exact sequence
$$0\rightarrow H \rightarrow G \rightarrow G/H \rightarrow 0$$ where
$H$ and $G$ are compactly generated and $H$ is a normal closed
subgroup of $G$. We found that actually we not need to assume that
$H$ and $G$ are compactly generated. We recall briefly the Milnor's
construction (see \cite{Mil}) and the Borel construction.

For a topological group $G$, the $n$-th join of $G$ is the space
$$G_n:=\{t_1g_1\oplus \cdots \oplus t_ng_n|\sum^n_{i=1}t_i=1,
t_i\geq 0, g_i\in G, \forall i=1,..., n\}/\sim$$ with the Milnor's
strong topology, and the equivalence relation $\sim$ identifies two
elements $t_1g_1\oplus \cdots \oplus t_ng_n\sim t'_1g'_1\oplus
\cdots \oplus t'_ng'_n$ if $t_i=t'_i$ for all $i=1,..., n$ and
$g_i=g'_i$ if $t_i\neq 0$.

We form the infinite join $E(G)$ of $G$ in the same manner, with the
restriction that all but a finite number of the $t_i$ should vanish.
Then $G$ acts freely on $E(G)$ and there exists a locally trivial
universal principal $G$-bundle
$$p: EG \rightarrow BG$$
where $BG=EG/G$. Furthermore, all the homotopy groups of $EG$
vanish, and from the homotopy sequence induced by the fibration $p$,
we see that $\pi_k(BG)=\pi_{k-1}(G)$ for $k\geq 1$.

If $\phi:G \rightarrow G'$ is a morphism between two topological
groups, define $E(\phi):E(G) \rightarrow E(G')$ by
$$E(\phi)(t_1g_1\oplus t_2g_2\oplus \cdots )=t_1\phi(g_1)\oplus
t_2\phi(g_2)\oplus \cdots$$ then $E(\phi)$ is an equivariant map and
induces a map $B(\phi):B(G) \rightarrow B(G')$.

\begin{definition}
Suppose that $G$ is a topological group acting on a topological
space $X$. The Borel construction is the orbit space  $B(X, G)=(X
\times EG)/G$ where  $EG$ is the  universal bundle of $G$. Then
there is a fibration
$$\xymatrix{ X \ar[r] & B(X, G) \ar[d]^{p} &\\
 & BG   \\}$$
where $p: B(X, G) \longrightarrow BG$ is the projection.
\end{definition}

\begin{definition}
We say that a map $\phi:A \rightarrow B$ between two topological
spaces is a weak homotopy equivalence if $\phi$ induces isomorphism
$\phi_*:\pi_n(A)\rightarrow \pi_n(B)$ for each $n\geq 0$. We say
that $A$, $B$ are weakly homotopy equivalent if they have the same
homotopy groups.
\end{definition}

\begin{proposition}
\begin{enumerate}
\item If $\phi:G \rightarrow G'$ is a topological group
homomorphism which is also a weak homotopy equivalence, then the
induced map $B(\phi):B(G) \rightarrow B(G')$ is also a weak homotopy
equivalence.

\item Let $G$ be a topological group and $\Omega BG$ be the loop space of $BG$
with based point $[e]\in BG$ where $e$ is the identity of $G$. Then
there is a map $\gamma:G \rightarrow \Omega BG$ which induces
isomorphisms of their homotopy groups.

\item Suppose that we have a commutative diagram of topological
group homomorphisms:
$$\xymatrix{ H_1 \ar[d]_{\psi_1} \ar[r]^{\phi_1} & G_1 \ar[d]^{\psi_2} &\\
H_2 \ar[r]^{\phi_2} & G_2}$$

Then $\psi_1, \psi_2$ induce a map $B(\psi_2, \psi_1): B(G_1, H_1)
\longrightarrow B(G_2, H_2)$. When $\psi_1, \psi_2$ are weak
homotopy equivalences, $B(\psi_2, \psi_1)$ is also a weak homotopy
equivalence.

\item If $\phi:G \rightarrow H$ is a morphism between topological
groups and $B(\phi)_*:B(G) \rightarrow B(H)$ is a weak homotopy
equivalence. If $G$ and $H$ have the homotopy type of a CW-complex,
then $\phi$ is a homotopy equivalence.

\item Suppose that $F:H_1 \times I \longrightarrow H_2$,
$F':G_1\times I \longrightarrow G_2$ are homotopies between
$F_0=\psi_1, F_1=\psi'_1$ and $F'_0=\psi_2, F'_1=\psi'_2$
respectively where each $F_t, F'_t$ are group homomorphisms. If the
following diagram commutes:
$$\xymatrix{ H_1 \ar[d]_{F_t} \ar[r]^{\phi_1} & G_1 \ar[d]^{F'_t} &\\
H_2 \ar[r]^{\phi_2} & G_2}$$ for each $t\in I$, then $B(\psi_2,
\psi_1)$ is homotopic to $B(\psi'_2, \psi'_1)$.
\end{enumerate}
\label{Borel construction}
\end{proposition}

\begin{proof}
\begin{enumerate}
\item Consider the long exact sequences on homotopy groups induced from
the two fibrations
$$\xymatrix{ \cdots \ar[r] & \pi_n(G) \ar[r] \ar[d]^{\phi_*} & \pi_n(E(G))
\ar[r] \ar[d]^{E(\phi)_*} & \pi_n(B(G)) \ar[r]
\ar[d]^{B(\phi)_*} & \cdots &\\
\cdots \ar[r] & \pi_n(G') \ar[r] & \pi_n(E(G')) \ar[r] &
\pi_n(B(G')) \ar[r] & \cdots &}$$ which shows that $B(\phi)_*$ is an
isomorphism between the homotopy groups of $B(G)$ and $B(G')$.

\item Let $PBG$ be the path space of $BG$ with based point $[e]\in
BG$. Define $\gamma: EG \rightarrow PBG$ by
$$\gamma(t_1g_1\oplus t_2g_2\oplus \cdots \oplus
t_ng_n):=[(1-t)e:tt_1g_1:tt_2g_2:\cdots :tt_ng_n]$$ where $0\leq
t\leq 1$. Then we have a commutative diagram of two fibrations:
$$\xymatrix{G \ar[r] \ar[d]^{\gamma} & EG \ar[r] \ar[d]^{\gamma} & BG \ar[d]^{=} \\
\Omega BG \ar[r] & PBG \ar[r] & BG}$$ where the second row is the
path fibration. Since $EG$ and $PBG$ are weakly contractible, the
map $\gamma$ induces isomorphisms between the homotopy groups of $G$
and $\Omega BG$.

\item
Define $B(\psi_2, \psi_1):G_1\times E(H_1)\rightarrow G_2\times
E(H_2)$ by $(x, y)\rightarrow (\psi_2(x), \psi_{1*}(y))$ which
induces a map $B(\psi_2, \psi_1):B(G_1, H_1) \rightarrow B(G_2,
H_2)$. Consider the long exact sequences
$$\xymatrix{ \cdots \ar[r] & \pi_n(G_1) \ar[r] \ar[d]^{\psi_{2*}} & \pi_n(B(G_1, H_1))
\ar[r] \ar[d]^{B(\psi_2, \psi_1)_*} & \pi_n(BH_1) \ar[r]
\ar[d]^{B\psi_{1*}} & \cdots &\\
\cdots \ar[r] & \pi_n(G_2) \ar[r] & \pi_n(B(G_2, H_2)) \ar[r] &
\pi_n(BH_2) \ar[r] & \cdots &}.$$ Since $\psi_{2*}$ and $B(\psi_1)$
are weak homotopy equivalence, $B(\psi_2, \psi_1)$ is also a weak
homotopy equivalence.

\item We have the following commutative diagram:
$$\xymatrix{ G \ar[r]^{\phi} \ar[d] & H \ar[d]\\
\Omega B(G) \ar[r] & \Omega B(H)}$$ where the vertical arrows are
weak homotopy equivalences. Hence the map $\phi$ induces
isomorphisms of homotopy groups between $G$ and $H$. If $G$ and $H$
have the homotopy type of a CW-complex, by the Whitehead theorem,
$\phi$ is a homotopy equivalence.

\item
By (ii), $F$ and $F'$ induce a map $B(F', F):B(G_1, H_1)\times I
\longrightarrow B(G_2, H_2)$ which is a homotopy between $B(F'_0,
F_0)=B(\psi_2, \psi_1)$ and $B(F'_1, F_1)=B(\psi'_2, \psi'_1)$.
\end{enumerate}
\end{proof}

The following technique is the main tool that we use to produce long
exact sequences of homotopy groups.
\begin{proposition}
If $H$ is a normal closed subgroup of a topological group $G$, then
there is a long exact sequence of homotopy groups:
$$\cdots \rightarrow \pi_n(H) \rightarrow \pi_n(G) \rightarrow
\pi_n(G/H) \rightarrow \pi_{n-1}(H) \rightarrow \cdots$$
\label{homotopy sequence}
\end{proposition}

\begin{proof}
From the fibrations
$$\xymatrix{H \ar[r] & E(G) \ar[d]&  E(G/H) \ar[r] & (E(G/H)\times
E(G))/G \ar[d]\\
& E(G)/H & & E(G)/G}$$ we see that $E(G)/H$ is a weak model of
$B(H)$ and $(E(G/H)\times E(G))/G$ is a weak model of $B(G)$. Since
$(E(G/H)\times E(G))/G=(E(G/H)\times (E(G)/H))/(G/H)$, by the Borel
construction, we have a fibration
$$\xymatrix{E(G)/H \ar[r] & (E(G/H)\times (E(G)/H))/(G/H) \ar[d]\\
& E(G/H)/(G/H)}$$ which induces the desire long exact sequence in
homotopy groups.
\end{proof}

\subsection{Relative Theory in Reduced Real Lawson Homology}
\begin{definition}
Suppose that $X, Y$ are real projective varieties and $Y$ is a
subvariety of $X$. Define the dimension $p$ relative reduced real
cycle group to be
$$R_p(X, Y)=\frac{R_p(X)}{R_p(Y)}.$$
\end{definition}

\begin{definition}
We say that a quasiprojective variety $U$ is real if there exist
real projective varieties $X$ and $Y$ , $Y\subset X$ such that
$U=X-Y$. We define the $p$-th reduced real cycle group of $U$ to be
$R_p(U)=R_p(X, Y)$.
\end{definition}

We need to show that the definition of reduced real cycle groups of
real quasiprojective varieties is independent of the choice of its
compactification.

\begin{lemma}
Suppose that $X, Y$ are real projective varieties and $Y$ is a
subvariety of $X$. Then $\frac{R_r(X)}{R_r(Y)}$ is isomorphic as a
topological group to $\frac{Z_r(X, Y)_{\R}}{Z_r(X, Y)^{av}}$.
\label{quotient isomorphic}
\end{lemma}

\begin{proof}
Let $Q_1, Q_2, q_1, q_2$ be the quotient maps and
$\psi(x+Z_r(X)^{av}+R_r(Y))=x+Z_r(Y)_{\R}+Z_r(X, Y)^{av}$ as
following:
$$\xymatrix{ & Z_r(X)_{\R} \ar[dl]_{q_1}   \ar[dr]^{q_2} & \\
R_r(X) \ar[d]_{Q_1} & & Z_r(X, Y)_{\R} \ar[d]^{Q_2} \\
\frac{R_r(X)}{R_r(Y)} \ar[rr]^{\psi}&  & \frac{Z_r(X,
Y)_{\R}}{Z_r(X, Y)^{av}}}$$ The diagram commutes and since those
quotient maps are open maps, $\psi$ is continuous. Easy to check
that $\psi$ is bijective. Since $\psi(U)=Q_2\circ q_2 \circ q_1^{-1}
\circ Q_1^{-1}(U)$ is open for an open set $U$, $\psi$ is an open
map hence the inverse of $\psi$ is also continuous.
\end{proof}

\begin{definition}
Suppose that $X, Y, X', Y'$ are real projective varieties, and
$Y\subset X$, $Y' \subset X'$. A regular map $f:X \rightarrow X'$ is
real if $f(\overline{x})=\overline{f(x)}$ for all $x\in X$. The pair
$(X, Y)$ is said to be relatively isomorphic to the pair $(X', Y')$
if there is a real regular map $f:X \longrightarrow X'$ such that
$f$ induces an isomorphism as real quasiprojective varieties between
$X-Y$ and $X'-Y'$, and we call $f$ a real relative isomorphism.
\end{definition}

In \cite[Theorem 4.3]{Lima}, Lima-Filho proved the following result.
\begin{theorem}
A relative isomorphism $f:(X, X') \rightarrow (Y, Y')$ induces an
isomorphism of topological groups:
$$f_*:\frac{Z_r(X)}{Z_r(X')} \rightarrow \frac{Z_r(Y)}{Z_r(Y')}$$
for all $r\geq 0$.
\end{theorem}

Hence when we have a real relative isomorphism, its restriction to
real cycle groups and averaged cycle groups gives us the following
isomorphisms.
\begin{proposition}
Suppose that $f:(X, Y) \longrightarrow (X', Y')$ is a real relative
isomorphism. Then $f$ induces a topological group isomorphism
between
\begin{enumerate}
    \item $Z_r(X, Y)_{\R}$ and $Z_r(X', Y')_{\R}$
    \item $Z_r(X, Y)^{av}$ and $Z_r(X', Y')^{av}$
    \item $\frac{Z_r(X, Y)_{\R}}{Z_r(X, Y)^{av}}$ and $\frac{Z_r(X',
    Y')_{\R}}{Z_r(X', Y')^{av}}$
\end{enumerate}
\label{Lima approach}
\end{proposition}

\begin{corollary}
Suppose that $f:(X, Y) \longrightarrow (X', Y')$ is a real relative
isomorphism, then $f$ induces a topological group isomorphism
between $\frac{R_r(X)}{R_r(Y)}$ and $\frac{R_r(X')}{R_r(Y')}$.
\end{corollary}

\begin{proof}
This follows directly from Lemma \ref{quotient isomorphic} and the
proposition above.
\end{proof}

Thus the reduced real cycle groups of a real quasiprojective variety
are well defined. Now we are able to define the reduced real Lawson
homology groups of a real quasiprojective variety.

\begin{definition}
Let $X, Y$ be two real projective varieties and $Y$ be a subvariety
of $X$. We define the relative reduced real Lawson homology groups
of $(X, Y)$ to be
$$RL_pH_n(X,Y)=\pi_{n-p}(R_p(X, Y)).$$

Suppose that $U$ is a real quasiprojective variety. Take a real
compactification $(X, Y)$ of $U$. We define the reduced real Lawson
homology groups of $U$ to be
$$RL_pH_n(U)=\pi_{n-p}(R_p(X, Y)).$$
\end{definition}

\begin{proposition}
Suppose that $X_3\subset X_2\subset X_1$ are real projective
varieties. There is a long exact sequence in reduced real Lawson
homology:
$$\cdots RL_pH_n(X_2, X_3)\longrightarrow RL_pH_n(X_1, X_3) \longrightarrow
RL_pH_n(X_1, X_2) \longrightarrow RL_pH_{n-1}(X_2,
X_3)\longrightarrow \cdots$$
\end{proposition}

\begin{proof}
We have a short exact sequence of topological groups
$$0\rightarrow \frac{R_p(X_2)}{R_p(X_3)} \rightarrow
\frac{R_p(X_1)}{R_p(X_3)} \rightarrow \frac{R_p(X_1)}{R_p(X_2)}
\rightarrow 0,$$ and by Proposition \ref{homotopy sequence}, we get
the long exact sequence of homotopy groups.
\end{proof}

\begin{corollary}
Suppose that $U, V$ are real projective varieties and $V\subset U$
is closed. Then there is a (localization) long exact sequence in
reduced real Lawson homology:
$$\cdots RL_pH_n(V) \longrightarrow RL_pH_n(U) \longrightarrow
RL_pH_n(U-V)\longrightarrow RL_pH_{n-1}(V)\longrightarrow \cdots$$
\label{localization}
\end{corollary}

\begin{proof}
Let $U=X_1-X_3$ where $X_1, X_3$ are real projective varieties and
take the closure $\overline{V}$ of $V$ in $X_1$. Let
$X_2=\overline{V} \cup X_3$. Then $X_2-X_3=V$, $X_1-X_2=U-V$ and
$X_3\subset X_2 \subset X_1$. By the Proposition above, we get the
long exact sequence.
\end{proof}

\subsection{Relative Theory in Reduced Real Morphic cohomology}
\begin{definition}
Suppose that $X, Y$ are real projective varieties and the dimension
of $Y$ is $n$. We define the codimension $t$ real cocycle group of
$X$ with values in $Y$ to be
$$R^t(Y)(X)=R_{n-t}(Y)(X),$$
and the codimension $t$  reduced real cocycle group of $X$ to be
$$R^t(X)=\frac{R_0(\P^t)(X)}{R_0(\P^{t-1})(X)}.$$
\end{definition}

From Proposition \ref{subspace}, we know that $R_0(\P^{t-1})(X)
\hookrightarrow R_0(\P^t)(X)$ as a closed subgroup, thus $R^t(X)$ is
a topological group.

To define relative reduced real morphic cohomology, we need to show
that $R^t(X)$ acts on $R^t(Y)$ where $Y$ is a real subvariety of
$X$.
\begin{proposition}
Suppose that $X$, $Y$ are real projective varieties and $Y$ is a
subvariety of $X$. Then there is a topological group homomorphism
$$\phi: R^t(X) \longrightarrow R^t(Y).$$
\label{relative action}
\end{proposition}

\begin{proof}
We have a commutative diagram:
$$\xymatrix{Z_0(\P^t)(X)^{av} \ar[r] \ar[d] & Z_0(\P^t)(X)_{\R} \ar[r]
\ar[d] & R_0(\P^t)(X) \\
Z_0(\P^t)(Y)^{av} \ar[r] & Z_0(\P^t)(Y)_{\R} \ar[r]& R_0(\P^t)(Y)}$$
where maps in the horizontal direction are inclusions and maps in
the vertical direction are restrictions. Hence they induce a map
from $R_0(\P^t)(X)$ to $R_0(\P^t)(Y)$.

From the following commutative diagram:
$$\xymatrix{R_0(\P^{t-1})(X) \ar[r] \ar[d] & R_0(\P^t)(X) \ar[r] \ar[d] &
R^t(X)\\
R_0(\P^{t-1})(Y) \ar[r]  & R_0(\P^t)(Y) \ar[r] & R^t(Y)},$$ we see
that there is a map from $R^t(X)$ to $R^t(Y)$.
\end{proof}

\begin{definition}
Suppose now $Y$ is a subvariety of $X$. Define the codimension $t$
relative reduced real cocycle group to be
$$R^t(X|Y)=B(R^t(Y), R^t(X))$$ by the action from Proposition
\ref{relative action}.
\end{definition}

\begin{definition}
The codimension $t$ reduced real morphic cohomology group of a real
projective variety $X$ is defined to be
$$RL^tH^k(X)=\pi_{t-k}(R^t(X))$$
Suppose now $X, Y$ are real projective varieties and $Y\subset X$.
We define the relative reduced real morphic cohomology group to be
$$RL^tH^k(X|Y)=\pi_{t-k}(R^t(X|Y)),$$
and the reduced real bivariant $Y$-valued morphic cohomology to be
$$RL^tH^k(X; Y)=\pi_{t-k}R^t(Y)(X).$$
\end{definition}

\begin{theorem}(Long exact sequences)
Long exact sequence in the reduced real morphic cohomology:
$$\cdots \longrightarrow RL^tH^k(Y) \longrightarrow RL^tH^k(X|Y)
\longrightarrow RL^tH^{k+1}(X) \longrightarrow RL^tH^{k+1}(Y)
\longrightarrow \cdots.$$
\end{theorem}

\begin{proof}
This follows from the homotopy sequence induced by the Borel
construction
$$\xymatrix{R^t(Y) \ar[r] & B(R^t(Y), R^t(X)) \ar[d]\\
& B(R^t(X))}$$
\end{proof}

\begin{remark}
It is much more difficult to generalize the reduced real morphic
cohomology to real quasiprojective varieties. The generalization of
morphic cohomology to quasiprojective varieties was done by
Friedlander in \cite{F2}. It is quite possible that his approach may
be adopted, but in this paper what we mainly concern are projective
varieties, so we do not attempt to do that now.
\end{remark}

\section{Functoriality}
There are analogous functorial properties between Lawson homology
and reduced real Lawson homology, morphic cohomology and reduced
real morphic cohomology. Similar results for bivariant morphic
cohomology can be found in \cite{FL1}.

\begin{proposition}
Suppose that $X, X', Y, Y'$ are real projective varieties and $f:X'
\longrightarrow X$, $g:Y\longrightarrow Y'$ are morphisms of real
projective varieties.

(a) The ``pullback of reduced real cocycles" determines a
homomorphism $$f^*:RL^sH^q(X; Y) \longrightarrow RL^sH^q(X'; Y),$$
$$f^*:RL^sH^q(X) \rightarrow RL^sH^q(X').$$

(b) The ``pushforward of reduced real cycles" determines a
homomorphism $$g_*:RL^sH^q(X; Y) \longrightarrow RL^{s-c}H^{q-c}(X;
Y'),$$ $$g_*:RL_rH_p(Y) \rightarrow RL_rH_p(Y')$$ where $c=dim Y-
dimY'$.

(c) Given morphisms $f_2: X' \longrightarrow X$,
$f_1:X''\longrightarrow X'$ and $g_1: Y\longrightarrow Y'$,
$g_2:Y'\longrightarrow Y''$ where $X, X', X'', Y, Y', Y''$ are all
real projective varieties, then
$$(f_2\circ f_1)^*=(f_1)^*\circ (f_2)^* \mbox{ and } (g_2\circ
g_1)_*=(g_2)_*\circ (g_1)_*.$$\label{pullback}
\end{proposition}

\begin{proof}
(a) The pullback $f^*: Z^s(Y)(X)\longrightarrow Z^s(Y)(X')$ is a
continuous map which is proved in \cite[Proposition 2.4]{FL1}. So we
only need to verify that $f^*$ maps real cycles to real cycles and
averaged cycles to averaged cycles. Since $f$ is real morphism, if
$\alpha\in Z^s(Y)(X)_{\R}$,
$f^*\alpha(\overline{x})=\alpha(f(\overline{x}))=\overline{\alpha(f(x))}$,
hence $f^*\alpha\in Z^s(Y)(X')_{\R}$. For $\beta\in Z^s(Y)(X)$,
$f^*(\beta+\overline{\beta})=f^*\beta+\overline{f^*\beta}\in
Z^s(Y)(X')^{av}$. So it induces a map from $R^s(Y)(X)$ to
$R^s(Y)(X')$, and therefore a group homomorphism in homotopy groups.
From the commutative diagram
$$\xymatrix{R_0(\P^{t-1})(X) \ar[r]^{f^*} \ar[d] & R_0(\P^{t-1})(X') \ar[d]\\
R_0(\P^t)(X) \ar[r] & R_0(\P^t)(X')}$$ we see that $f^*$ induces a
map $f^*:R^t(X) \rightarrow R^t(X')$.

(b) Let $n$ be the dimension of $Y$. The pushforward $g_*$
determines a continuous algebraic monoid homomorphism
$$g_*:\mathscr{C}_{n-s}(Y) \rightarrow \mathscr{C}_{n-s}(Y'),$$
and by composition, a continuous homomorphism
$$\mathcal{M}or(X, \mathscr{C}_{n-s}(Y)) \rightarrow
\mathcal{M}or(X, \mathscr{C}_{n-s}(Y')).$$ Passing to group
completion we get a map $g_*:Z^s(Y)(X) \rightarrow Z^{s-c}(Y')(X)$.
It is easy to see that $g_*$ maps real cycles to real cycles and
averaged cycles to averaged cycles. Hence we get a map
$g_*:R^s(Y)(X) \rightarrow R^{s-c}(Y')(X)$. When $X$ is a point, it
reduces to the case $g_*:R_r(Y) \rightarrow R_r(Y')$.

(c) Since $(f_2\circ f_1)^*=f^*_1\circ f^*_2$ for $Y$-valued cocycle
groups as proved in \cite[Proposition 2.4]{FL1}, from (a), this
relation passes to the quotient. Thus $(f_2\circ f_1)^*=f^*_1\circ
f^*_2$ in reduced $Y$-valued real cycle groups and reduced real
cycle groups. A similar argument works for $g_1, g_2$.
\end{proof}

\begin{proposition}
Let $f:X\longrightarrow X'$, $g:Y'\longrightarrow Y$ be morphisms of
real projective varieties.

(a) If $f$ has equidimensional fibers, then there are Gysin
homomorphisms $$f_!:RL^sH^q(X; Y) \longrightarrow
RL^{s-c}H^{q-c}(X';Y),$$ $$f_!:RL^sH^q(X) \rightarrow
RL^{s-c}H^{q-c}(X')$$ for all $s, q$ with $s\geq q \geq c$ where
$c=dim(X)-dim(X')$.

(b) If $g$ is flat, then for all $s, q$ with $s\geq q$, there are
Gysin homomorphisms $$g^!:RL^sH^q(X; Y)\longrightarrow RL^sH^q(X;
Y'),$$ $$g^!:RL_rH_p(Y) \rightarrow RL_{r+e}H_{p+e}(Y')$$ where
$e=dim(Y')-dim(Y)$.

(c). If $f_1: X \longrightarrow X'$ and $f_2:X' \longrightarrow X''$
are as in part (a), or if $g_2:Y'\longrightarrow Y$ and
$g_1:Y''\longrightarrow Y'$ are as in part (b), then
$$(f_2\circ f_1)_!=(f_2)_!\circ(f_1)_! \mbox{ and } (g_2 \circ
g_1)^!=(g_1)^!\circ (g_2)^!.$$
\end{proposition}

\begin{proof}
(a) Let $f_!$ be the composition of the following maps:
$$Z_{n-s}(Y)(X)\rightarrow Z_{n+m-s}(X\times Y) \overset{(1\times f)_*}{\longrightarrow} Z_{n-s+c}(X\times Y)(X')
\overset{p_*}{\longrightarrow} Z_{n-s+c}(Y)(X')$$ where $1\times
f:X\times Y \rightarrow X\times Y \times X'$ is the map defined by
sending $(x, y)$ to $(x, y, f(x))$ and $p:X\times Y \rightarrow Y$
is the projection. Then $f_!$ is a real map and hence induces a map
$f_!:R^s(Y)(X) \rightarrow R^{s-c}(Y)(X')$. By the Lawson suspension
theorem \ref{Lawson suspension}, we have a homotopy equivalence
between $R_c(\P^s)(X)$ and $R_0(\P^{s-c})(X)$. From the following
commutative diagram
$$\xymatrix{R_0(\P^{s-1})(X) \ar[r]^{f_!} \ar[d] & R_c(\P^{s-1})(X') \ar[d] \ar[r]^{\cong} & R_0(\P^{s-1-c})(X') \ar[d] \\
R_0(\P^s)(X) \ar[r]^{f_!} & R_c(\P^s)(X') \ar[r]^{\cong} &
R_0(\P^{s-c})(X')}$$ we see that $f_!$ induces a map $f_!:R^s(X)
\rightarrow R^{s-c}(X')$.

(b) Since $g$ is flat, we have a map $g^!:\mathscr{C}_{n-s}(Y)
\rightarrow \mathscr{C}_{n-s+e}(Y')$ of pullback cycle. Then by
composition, we get a map $g^!:\mathscr{M}or(X,
\mathscr{C}_{n-s}(Y)) \rightarrow \mathscr{M}or(X,
\mathscr{C}_{n-s+e}(Y'))$ which induces a map after group completion
$g^!:Z^s(Y)(X) \rightarrow Z^s(Y')(X)$. It is obvious that $g^!$ is
real, hence it induces a map $g^!:R^s(Y)(X) \rightarrow R^s(Y')(X)$.
Take $X$ to be a point, we get $g^!:R_r(Y) \rightarrow R_{r+e}(Y')$.
Taking homotopy groups, we complete the proof.
\end{proof}

\section{Fundamental Properties}
\subsection{Homotopy Invariance}
\begin{proposition}
Let $X, X'$ and $Y$ be real projective varieties. Suppose that
$F:I\times X \longrightarrow X'$ is a continuous map where $I$ is
the unit interval and for each $t\in I$ , $F_t$ is a real algebraic
morphism. If $F_0(x)=f(x), F_1(x)=g(x)$ for all $x\in X$, then
$f^*=g^*:RL^*H^*(X'; Y) \longrightarrow RL^*H^*(X; Y)$ and they
induce the same map $f^*=g^*:RL^*H^*(X') \longrightarrow
RL^*H^*(X)$.
\end{proposition}

\begin{proof}
Consider the map $$F^*:Z_r(Y)(X') \times \C \longrightarrow
Z_r(Y)(X)$$ which sends $(\sigma, t)$ to $F^*_t\sigma$. This map is
continuous and is a homotopy between $f^*$ and $g^*$. Since each
$F_t$ is a real morphism, by Proposition \ref{pullback}, $F^*_t$
passes to the quotient and thus $F$ induces a map $F^*:
R_r(Y)(X')\times I \longrightarrow R_r(Y)(X)$ which is a homotopy
between $F^*_0=f^*$ and $F^*_1=g^*$ so they induce the same map from
$RL^*H^*(X'; Y)$ to $RL^*H^*(X; Y)$. Furthermore, $F^*$ induces a
homotopy $R_0(\P^{t-1})(X')\times I \longrightarrow
R_0(\P^{t-1})(X)$ and a homotopy $R_0(\P^t)(X')\times I
\longrightarrow R_0(\P^t)(X)$ between the two maps induced by $f$
and $g$. Hence their induced maps from $R^t(X')$ to $R^t(X)$ are
homotopic.
\end{proof}

\subsection{The Splitting Principle}
\begin{proposition}
The real projective space $\RP^d$ is homeomorphic to
$SP^d(\P^1)_{\R}$ where $SP^d(\P^1)_{\R}$ is the subset of
$SP^d(\P^1)$ consisting of zero cycles invariant under the
conjugation.
\end{proposition}

\begin{proof}
A proof of $\P^d\cong SP^d(\P^1)$ can be found in [12]. The main
idea is sketched in the following. For $p=\{p_1,..., p_d\}\in
SP^d(\P^1)$, let $p_i=[-b_i:a_i]\in \P^1$ and we associate $p$ to a
homogeneous polynomial
$$P(x, y)=\prod^d_{i=1}(a_ix+b_iy)=\sum^d_{k=0}c_kx^ky^{d-k}$$
where $$c_k=\sum_{|I|=k}a_Ib_{I'}$$ and the sum is taken over all
multi-indices $I=\{0\leq i_1 < \cdots < i_k \leq d\}$ of length
$|I|=k$ and $I'$ is the complementary multi-index with $|I'|=d-k$.
The map $\psi :SP^d(\P^1) \longrightarrow \P^d$ which maps $p$ to
the point $[c_0:\cdots: c_d]\in \P^d$ is an isomorphism. The
conjugation on $SP^d(\P^1)$ is given by the the conjugation of
$\P^d$ through this isomorphism, i.e., $\overline{p}$ is defined to
be $\psi^{-1}\overline{\psi(p)}$. Hence $\psi$ is real. This
conjugation is same as the conjugation induced by the conjugation of
$\P^1$ on $SP^d(\P^1)$, i.e., $\bar{p}=\{\bar{p_1},\cdots,
\bar{p_d}\}$. Since $\psi$ is a real homeomorphism, its restriction
to $SP^d(\P^1)_{\R}$ gives a homeomorphism between $SP^d(\P^1)_{\R}$
and $\RP^d$.
\end{proof}

From section 2.10 of \cite{FL1}, for $1\leq t \leq s$, there is a
monoid homomorphism
$$SP(\P^s) \overset{p}{\longrightarrow} SP(\P^t)$$
which is induced by the morphism
$$p: \P^s=SP^s(\P^1) \longrightarrow SP^{\binom{s}{t}}SP^t(\P^1)$$
defined by
$$\{x_1,..., x_s\} \longmapsto \underset{|I|=t}{\sum}\{x_{i_1},..., x_{i_t}\}.$$

Let $\P^0=\{x_0\}$. We have an embedding $SP^t(\P^1) \hookrightarrow
SP^s(\P^1)$ given by
$$\{x_1,..., x_t\} \mapsto \{x_1,..., x_t, x_0,..., x_0\}$$

We extend this morphism additively to a monoid homomorphism
$$p: SP(\P^s) \rightarrow SP(\P^t)$$
and extend this morphism further to include the case $t=0$ by
setting
$$p(\sum n_{\alpha} x_{\alpha})=(\sum n_{\alpha})x_0$$ for
$x_{\alpha}\in \P^s$ and $n_{\alpha} \in \Z^+$. It is then easy to
see that $p$ is real.

The map $p$ induces a homomorphism
$$p: \mathscr{C}_0(\P^s)(X) \rightarrow \mathscr{C}_0(\P^t)(X)$$
which passes to the group completion
$$p: Z_0(\P^s)(X) \rightarrow Z_0(\P^t)(X)$$
Since $p$ is real, the map also induces two maps
$$p: Z_0(\P^s)(X)_{\R} \rightarrow Z_0(\P^t)(X)_{\R}$$
and
$$p: Z_0(\P^s)(X)^{av} \rightarrow Z_0(\P^t)(X)^{av}$$ hence
a map
$$p': R_0(\P^s)(X) \rightarrow R_0(\P^t)(X)$$ for $0\leq t \leq s$.

The inclusion map $i: \P^t \longrightarrow \P^s$ induces an
inclusion map $i': R_0(\P^t)(X) \longrightarrow R_0(\P^s)(X)$ by
Proposition \ref{subspace}.

\begin{lemma}
The composition $\varphi'=p'\circ i' :R_0(\P^t)(X) \longrightarrow
R_0(\P^t)(X)$ is of the form
$$\varphi'=Id+\psi'$$ where image $\psi'\subset R_0(\P^{t-1})(X)$.
\end{lemma}

\begin{proof}
Let $\varphi=p \circ i$ and $\psi=\varphi-Id$. Since $\varphi$ is
real, $\psi$ is real. By \cite[Lemma 2.11]{FL1}, $p\circ i
(f)=(Id+\psi)(f)$ for $f\in Z_0(\P^t)(X)$ where the image of $\psi$
is contained in $Z_0(\P^{t-1})(X)$. By the definition of $p'$ and
$i'$, we have $p'\circ i'(f+Z_0(\P^t)(X)^{av})=p\circ i
(f)+Z_0(\P^t)(X)^{av}=f+\psi(f)+Z_0(\P^t)(X)^{av}$ for $f\in
Z_0(\P^t)(X)_{\R}$. The map $\psi$ is real, thus
$\psi(Z_0(\P^t)(X)_{\R})$ is contained in $Z_0(\P^{t-1})(X)_{\R}$.
And it induces a map $\psi': R_0(\P^t)(X) \rightarrow R_0(\P^t)(X)$
such that $p'\circ i'=(Id+\psi')$ and the image of $\psi'$ is
contained in $R_0(\P^{t-1})(X)$.
\end{proof}

\begin{theorem}(Splitting Principle)
For a real projective variety $X$, we have a homotopy equivalence:
$$R_0(\P^s)(X) \cong R^0(X)\times R^1(X) \times \cdots \times
R^s(X).$$
\end{theorem}

\begin{proof}
Let $M^t=R_0(\P^t)(X)$, $Q^t=\frac{M^t}{M^{t-1}}$. Let $p^{k, t}:
M^k \rightarrow M^t$ be the map induced by $p$, and $q^t:M^t
\rightarrow Q^t$ be the quotient map. Let $q^{k, t}:M^k \rightarrow
Q^t$ be the map defined by $q^{k, t}=q^t\circ p^{k, t}$, and
$\xi^{k, t}:M^t \hookrightarrow M^k \rightarrow Q^0\times \cdots
Q^t$ defined by $\xi^k=(q^{k, 0}\circ i)\times \cdots \times (q^{k,
t}\circ i)$. We use induction to show that the map
$$M^t \overset{\xi^{k, t}}{\longrightarrow}
Q^0\times \cdots Q^t$$ is a homotopy equivalence. When $t=0$, it is
trivially true. We assume that it is true for $t-1$. We have a
commutative diagram:
$$\xymatrix{M^{t-1} \ar[r] \ar[d]^{\xi^{k, t-1}} & M^t \ar[r]^{q^t}
\ar[d]^{\xi^{k, t}} & Q^t \ar[d]^{=}\\
Q^0\times \cdots \times Q^{t-1} \ar[r] & Q^0\times \cdots \times
Q^{t-1}\times Q^t \ar[r] & Q^t}$$ By Proposition \ref{Borel
construction}, we see that $\xi^{k, t}$ induces an isomorphism of
between homotopy groups of $M^t$ and $Q^0\times \cdots \times Q^t$.
In Appendix, we show that $M^t$ and $Q^0\times \cdots \times Q^t$
have the homotopy type of a CW-complex, hence by Whitehead Theorem,
$\xi^{k, t}$ is a homotopy equivalence.
\end{proof}

\subsection{The Lawson Suspension Theorem}
Let $\P^0=\{x_\infty\}$. The suspension $\susp X\subset \P^{n+1}$ of
a projective variety $X$ is the complex cone over $X$, or
equivalently, the Thom space of the hyperplane bundle
$\mathcal{O}(1)$ of $\P^n$ restricted to $X$. A point in $\susp X$
can be written as $[t:x]$ or $[1:0:...:0]$ where $t\in \C, x\in X$.
We consider $X$ as a subvariety of $\susp X$ by identifying $X$ with
the zero section of $\susp X$. If $X$ is a real projective variety,
$\susp X$ is also a real projective variety. The following
Proposition is straight forward.

\begin{proposition}
\begin{enumerate}
\item
Suppose that $\varphi: X\longrightarrow Y$ is a morphism of
projective varieties where $X\subset \P^n, Y\subset \P^m$. Then the
following diagram commutes:
$$\xymatrix{X \ar[r]^{\varphi} \ar[d] & Y \ar[d]\\
\susp X \ar[r]^{\widetilde{\susp}\varphi} & \susp Y}$$ where
$(\widetilde{\susp}\varphi)([x:t])=[\varphi(x):t]$ for all $t\in \C,
x\in X$ and $(\widetilde{\susp}\varphi)([0:...:0:1])=[0:...:0:1]$
where the vertical arrows are inclusions.

\item
If $\varphi: X \rightarrow Y$ is a morphism between real projective
varieties, then $\widetilde{\susp}\varphi$ is real.
\end{enumerate}
\end{proposition}

The above commutativity induces a commutative diagram of cycle
groups.
\begin{corollary}
\begin{enumerate}
\item Suppose that $\varphi: X\longrightarrow Y$ is a morphism of
projective varieties. Then the following diagram commutes:
$$\xymatrix{Z_r(X) \ar[r]^{\varphi_*} \ar[d] & Z_r(Y) \ar[d]\\
Z_r(\susp X) \ar[r]^{(\widetilde{\susp}\varphi)_*} & Z_r(\susp
Y)}.$$

\item Suppose that $\varphi: X\longrightarrow Y$ is a morphism of
real projective varieties. Then the following diagram commutes:
$$\xymatrix{R_r(X) \ar[r]^{\varphi_*} \ar[d] & R_r(Y) \ar[d]\\
R_r(\susp X) \ar[r]^{(\widetilde{\susp}\varphi)_*} & R_r(\susp Y)}$$
\end{enumerate}
\end{corollary}

For $f\in \mathscr{M}or(X, \mathscr{C}_r(Y))$, define the suspension
map $\susp_* : \mathscr{M}(X, \mathscr{C}_r(Y)) \longrightarrow
\mathscr{M}or(X, \mathscr{C}_{r+1}(\susp Y))$ by $(\susp_*
f)(x)=\susp f(x)$ for $x\in X$, the point-wise suspension. Passing
to the group completion, we get a map
$$\susp_*: Z_r(Y)(X) \rightarrow Z_{r+1}(\susp Y)(X)$$
The Lawson suspension theorem for bivariant morphic cohomology is
proved in \cite[Theorem 3.3]{FL1}. We state it in the following.

\begin{theorem}
For two projective varieties $X, Y$, the suspension map $\susp_{*} :
Z_r(Y)(X) \longrightarrow Z_{r+1}(\susp Y)(X)$ is a homotopy
equivalence.
\end{theorem}

The Lawson suspension theorem for real cycle groups and averaged
groups of projective spaces can be found in \cite{Lam}. The
equivariant version of the Lawson suspension theorem can be found in
\cite{LLM2}. We generalize all previous results to bivariant reduced
real morphic cohomology.

Suppose that $V$ is a projective variety. Two subvarieties of $V$
meet properly in $V$ if the codimension of each component of their
intersection is the sum of the codimension in $V$ of that two
subvarieties. When we say that a cycle $c_1$ meets another cycle
$c_2$ properly we mean that each component of $c_1$ and $c_2$ meets
properly.

\begin{theorem}(Lawson suspension theorem)
Suppose that $X, Y$ are real projective varieties and $Y\subset
\P^N$. Then
\begin{enumerate}
\item $\susp_*: Z_r(Y)(X)_{\R} \longrightarrow Z_{r+1}(\susp
Y)(X)_{\R}$ is a homotopy equivalence.

\item $\susp_*: Z_r(Y)(X)^{av} \longrightarrow Z_{r+1}(\susp
Y)(X)^{av}$ is a homotopy equivalence.

\item $\susp_* : R_r(Y)(X)\longrightarrow R_{r+1}(\susp Y)(X)$ is
a homotopy equivalence.
\end{enumerate}
\label{Lawson suspension}
\end{theorem}

\begin{proof}
\item Step 1:

Let $$T_{r+1}(\susp Y)(X)=\{ f\in Z_{r+1}(\susp Y)(X)| f(x) \mbox{
meets } x \times Y \mbox{ properly for all } x\in X\},$$
$$T_{r+1}(\susp Y)(X)_{\R}=T_{r+1}(\susp Y)(X)\cap Z_{r+1}(\susp
Y)(X)_{\R},$$
$$T_{r+1}(\susp Y)(X)^{av}=T_{r+1}(\susp Y)(X)\cap Z_{r+1}(\susp
Y)(X)^{av}.$$

Following Proposition 3.2 in \cite{F1}, we define $\Lambda \subset
\P^{N+1}\times \P^1 \times \P^{N+1}$ to be the graph of the rational
map $\P^{N+1}\times \P^1 \longrightarrow \P^{N+1}$ whose restriction
to $\P^{N+1}\times t$ for $t\in \A^1-\{0\}$ is the linear
automorphism $\Theta_t:\P^{N+1}\longrightarrow \P^{N+1}$ sending
$[z_0:\cdots :z_N: z_{N+1}]$ to $[z_0:\cdots :z_N:
\frac{1}{t}z_{N+1}]$. More explicitly, $\Lambda$ is the closed
subvariety given by the homogeneous equations:
$$\left\{%
\begin{array}{ll}
    X_iY_j-X_jY_i=0 \\
    TX_{N+1}Y_j-SX_jY_{N+1}=0, & \hbox{ for } 0\leq i, j\leq N. \\
\end{array}%
\right.$$ where $[X_0:\cdots:X_{N+1}]\in \P^{N+1}, [S:T]\in \P^1$,
and $[Y_0:\cdots: Y_{N+1}]\in \P^{N+1}$.

For $f\in T_{r+1}(\susp Y)(X)$, $t\in \R-\{0\}\subset \P^1$, define
$$\phi_t(f)=Pr_{1, 4*}[(X\times \Lambda)\bullet (f\times t \times
\P^{N+1})]$$ where $Pr_{1, 4}:X\times \P^{N+1}\times \P^1\times
\P^{N+1} \longrightarrow X\times \P^{N+1}$ is the projection and it
is proved in \cite[Proposition 3.2]{F1}, $\phi_t(f)(x)$ meets $x
\times Y$ properly in $x\times \susp Y$ for all $x\in X$. For $t\in
\R$, since $\phi_t$ is real, $\phi_t$ maps $T_{r+1}(\susp
Y)(X)_{\R}$ to $T_{r+1}(\susp Y)(X)_{\R}$ and $T_{r+1}(\susp
Y)(X)^{av}$ to $T_{r+1}(\susp Y)(X)^{av}$.

When $t=0$, check the equations defining $\Lambda$, we see that
$\phi_0(f)=\susp(f\bullet (X\times Y))\in \susp_*Z_r(Y)(X)_{\R}$.
When $t=1$, $\phi_1$ is the identity map. If $f\in
\susp_*Z_r(Y)(X)$, $\phi_t(f)=f$ for $t\in \R$.

Therefore, $\phi_0: T_{r+1}(\susp Y)(X)_{\R} \longrightarrow
\susp_*Z_r(Y)(X)_{\R}$ and $\phi_0: T_{r+1}(\susp Y)(X)^{av}
\longrightarrow \susp_*Z_r(Y)(X)^{av}$ are strong deformation
retractions.

\item Step 2:

We use the notation $K_{r+1, e}(\susp Y)(X)$ as in section 2. Let
$$T_{r+1, e}(\susp Y)(X)=K_{r+1, e}(\susp Y)(X)\cap T_{r+1}(\susp
Y)(X),$$

$$T_{r+1, e}(\susp Y)(X)_{\R}=T_{r+1}(\susp Y)(X)_{\R}\cap T_{r+1,
e}(\susp Y)(X),$$ and $$T_{r+1, e}(\susp Y)(X)^{av}= T_{r+1}(\susp
Y)(X)^{av}\cap T_{r+1, e}(\susp Y)(X).$$

The second step is to show that the inclusion map $T_{r+1}(\susp
Y)(X)_{\R}\hookrightarrow Z_{r+1}(\susp Y)(X)_{\R}$ and the
inclusion map $T_{r+1}(\susp Y)(X)^{av} \hookrightarrow
Z_{r+1}(\susp Y)(X)^{av}$ induce homotopy equivalences. We follow
the proof of \cite[Theorem 3.3]{FL1}. Take two real points
$t_0=[0:\cdots:1], t_1=[0:\cdots:1:1]\in \P^{N+2}-\P^{N+1}$.
Consider $\P^{N+2}$ as the algebraic suspension of $\P^{N+1}$ with
respect to the point $t_0$ and let $pr:
\P^{N+2}-\{t_1\}\longrightarrow \P^{N+1}$ denote the linear
projection away from $t_1$. Then consider the partially defined
function
$$\Psi_e: K_{r+1, d}(\susp Y)(X)\times Div_e \longrightarrow \
K_{r+1, de}(\susp Y)(X)$$ given by
$$\Psi_e(f, D)=(Id\times pr)((\susp f)\bullet (X\times D))$$ where
$Div_e$ is the set of effective real divisors of degree $e$ on
$\P^{N+2}$.

Let
$$\bigtriangleup(f(x))=\{D\in Div_e: (\susp f(x) \bullet D) \mbox{ does not meet }
\{x\}\times H \mbox{ properly}\}$$ where $H\subset \P^{N+2}$ is the
hyperplane containing $\P^N\cup \{t_1\}$ and $f \in K_{r+1, e}(\susp
Y)(X)$.

By \cite[Lemma 5.11]{Lawson1}, one has
$$codim_{\R}(\bigtriangleup(f(x)))\geq \binom{p+e+1}{e}$$
where $p=dim f(x)$. In particular, $codim_{\R}(\bigtriangleup
(f(x))) \longrightarrow \infty$ as $e\longrightarrow \infty$.

We now choose $e(d)$ so that
$codim_{\R}(\bigtriangleup(f(x)))>dimX+1$ for all $e>e(d)$. Then
setting
$$\bigtriangleup (f)=\cup_{x \in X} \bigtriangleup(f(x)),$$
we have $codim_{\R}(\bigtriangleup (f))>1$ for all $e>e(d)$.
Consequently, for each $e>e(d)$, there must exist a line $L_e\subset
Div_e$ containing $e\cdot \P^{N+1}$ such that
$$(L_e-e\cdot \P^{N+1})\cap \bigtriangleup (f)=\emptyset.$$
It follows immediately that $\Psi_e$ restricted to $K_{r+1, e}(\susp
Y)(X)\times (L_e-e\cdot \P^{N+1})$ has images in $T_{r+1, de}(\susp
Y)(X)$ and this map is homotopic to the map given by multiplying by
$e$. Since $\Psi_e$ is a real map, by restricting $\Psi_e$, we get
two maps:
$$\Psi_e: K_{r+1, d}(\susp Y)(X)_{\R}\times (L_e-e\cdot \P^{N+1}) \longrightarrow \
T_{r+1, de}(\susp Y)(X)_{\R},$$
$$\Psi_e: K_{r+1, d}(\susp Y)(X)^{av} \times (L_e-e\cdot \P^{N+1}) \longrightarrow \
T_{r+1, de}(\susp Y)(X)^{av}$$ which are both homotopic to the map
given by multiplying by $e$. It then follows as in \cite[Theorem
4.2]{Lawson1} that the inclusions induce weak homotopy equivalences.

Consider the following commutative diagram:
$$\xymatrix{Z_r(Y)(X)^{av} \ar[r] \ar[d]^{\susp_*} & Z_r(Y)(X)_{\R}
\ar[d]^{\susp_*} \ar[r] & R_r(Y)(X) \ar[d]^{\susp_*}\\
Z_{r+1}(\susp Y)(X)^{av} \ar[r] & Z_{r+1}(\susp Y)(X)_{\R} \ar[r] &
R_{r+1}(Y)(X)}$$ where the first two vertical arrows are homotopy
equivalences induced by the suspension map. By Proposition
\ref{Borel construction}, there is a weak homotopy equivalence
induced by $\susp $ from $R_r(Y)(X)$ to $R_{r+1}(\susp Y)(X)$. In
Appendix we show that $Z_r(Y)(X)_{\R}, Z_r(Y)(X)^{av}$, and
$R_r(Y)(X)$ have the homotopy type of a CW-complex. Hence by the
Whitehead theorem, $\susp_*$ is a homotopy equivalence.
\end{proof}

Take $X$ to be a real point, we have the Lawson Suspension Theorem
for reduced real Lawson homology.

\begin{corollary}
$R_t(X)$ is homotopy equivalent to $R_{t+1}(\susp X)$ for $t\geq 0$.
\end{corollary}

\begin{corollary}
$R_t(\P^n)$ is homotopy equivalent to $K(\Z_2, 0)\times K(\Z_2, 1)
\times \cdots \times K(\Z_2, n-t)$.
\end{corollary}

\begin{proof}
By the Lawson suspension theorem \ref{Lawson suspension}, we have
$R_t(\P^n)\cong R_0(\P^{n-t})$. Since
$R_0(\P^{n-t})=\frac{Z_0(\RP^{n-t})}{2Z_0(\RP^{n-t})}$ where
$Z_0(\RP^{n-t})$ is the free abelian group generated by points of
$\RP^{n-t}$. By the Dold-Thom theorem,
$\frac{Z_0(\RP^{n-t})}{2Z_0(\RP^{n-t})}$ is homotopy equivalent to
$K(\Z_2, 0)\times K(\Z_2, 1) \times \cdots \times K(\Z_2, n-t)$.
\end{proof}

The following is the Lawson suspension theorem for reduced real
cocycle groups.
\begin{corollary}
$R^t(X)$ is homotopy equivalent to $\frac{R_p(\P^{p+t+1})(X)}{
R_p(\P^{p+t})(X)}$.
\end{corollary}

\begin{proof}
We have a commutative diagram where the vertical columns are
homotopy equivalences from the suspension theorem:
$$\xymatrix{R_0(\P^{t-1})(X) \ar[r] \ar[d] & R_0(\P^t)(X)
\ar[d]\\
R_p(\P^{p+t})(X) \ar[r] & R_p(\P^{p+t+1}(X))}$$ and then by
Proposition \ref{Borel construction}, $R^t(X)$ is homotopy
equivalent to $\frac{R_p(\P^{p+t+1})}{R_p(\P^{p+t}(X))}$.
\end{proof}

\begin{corollary}
The group $R_r(\C^n)$ is homotopy equivalent to $K(\Z_2, n-r)$ for
$r\leq n$, $n\in \N$. \label{C^n}
\end{corollary}

\begin{proof}
We have a commutative diagram:
$$\xymatrix{ R_r(\P^{n-1}) \ar[r] \ar[d]^{\susp_*} & R_r(\P^n) \ar[r] \ar[d]^{\susp_*} &
R_r(\C^n) \ar[d]\\
R_0(\P^{n-r-1}) \ar[r] & R_0(\P^{n-r}) \ar[r] & R_0(\C^{n-r})}$$ By
Proposition \ref{Borel construction} and the Dold-Thom theorem, we
see that $R_r(\C^n)$ is homotopy equivalent to $K(\Z_2, n-r)$.
\end{proof}

\subsection{Homotopy Property}
\begin{theorem}
Let $X$ be a real quasiprojective variety and $\pi: E \rightarrow X$
be a real algebraic vector bundle over $X$ of rank $k$, i.e., $E$ is
a real quasiprojective variety and $\pi: E \rightarrow X$ is a
complex algebraic vector bundle or rank $k$ with $\pi$ real. Then
$\pi^*: R_r(X) \longrightarrow R_{r+k}(E)$ induces a homotopy
equivalence.
\end{theorem}

\begin{proof}
If the dimension of $X$ is 0, then $E=\C^k$, this follows from
Corollary \ref{C^n}. Assume that the result is true for any real
quasiprojective variety of dimension $\leq m$. Let the dimension of
$X$ be $m+1$. Take a real quasiprojective variety $Y$ of dimension
$\leq m$ which is closed in $X$ such that $E$ restricted to $U=X-Y$
is trivial. We have the following commutative diagram of topological
groups:
$$\xymatrix{R_r(Y) \ar[r] \ar[d] & R_r(X) \ar[r] \ar[d] & R_r(U) \ar[d]\\
R_{r+k}(E|_Y) \ar[r] & R_{r+k}(E) \ar[r] & R_{r+k}(U\times \C^k)}$$

By the induction hypothesis, $\pi^*: R_r(Y)\longrightarrow
R_{r+k}(E|_Y)$ is a homotopy equivalence. Hence by Proposition
\ref{Borel construction}, the result follows from the case of
trivial bundles. By writing $\C^k=\C^{k-1}\times \C$, we see that it
is sufficient to prove for the case of trivial line bundles.

Let $\overline{X}$ be a real projective closure of $X$. Take real
projective varieties $Y, Y'\subset \overline{X}$ such that
$X=\overline{X}-Y$, $U=\overline{X}-Y'$ and the hyperplane line
bundle $\mathcal{O}_{\overline{X}}(1)|_U=U\times \C$. We have the
following commutative diagram of topological groups:
$$\xymatrix{R_r(Y\cup Y') \ar[d] \ar[r] & R_r(\overline{X}) \ar[r] \ar[d] & R_r(X\cap U)
\ar[d]\\
 R_{r+1}(\mathcal{O}_{Y\cup Y'}(1)) \ar[r] &
 R_{r+1}(\mathcal{O}_{\overline{X}}(1)) \ar[r] &
 R_{r+1}((X\cap U)\times \C)}$$ We note that for a real projective
variety $W$,
$$R_{r+1}(\mathcal{O}_W(1))=\frac{R_{r+1}(\overline{\mathcal{O}_{\overline{W}}(1)})}{R_{r+1}(\infty)}
=R_{r+1}(\susp W)\cong R_r(W)$$ Hence the first two vertical arrows
are homotopy equivalences. It follows that
$$\pi^*:R_r(X\cap U)\longrightarrow R_{r+1}((X\cap U)\times \C)$$ is a homotopy
equivalence. Since $X=(X\cap U)\cup (X\cap Y')$, we have the
following commutative diagram of topological groups:
$$\xymatrix{R_r(X\cap Y') \ar[d] \ar[r] & R_r(X) \ar[r] \ar[d] & R_r(X\cap U)\ar[d]\\
R_{r+1}((X\cap Y')\times \C) \ar[r] & R_{r+1}(X\times \C) \ar[r] &
R_{r+1}((X\cap U)\times \C)}$$ Since the dimension of $X\cap Y'$ is
less than the dimension of $X$, by induction on dimension, the left
vertical arrow is a homotopy equivalence. Therefore, by Proposition
\ref{Borel construction}, $\pi^*: R_r(X)\longrightarrow
R_{r+1}(X\times \C)$ is a homotopy equivalence.
\end{proof}

\section{Duality Theorem}
The main ingredient of the proof of the duality theorem between
morphic cohomology and Lawson homology is the moving lemma proved by
Friedlander and Lawson in \cite{FL2}. Here we observe that their
moving lemma preserves reality, and hence we have a moving lemma for
real cycle groups and averaged cycle groups. Since we are not going
to use the full power of the Friedlander-Lawson moving lemma, we
state a simple form which fits our need.

\begin{theorem}(Real Friedlander-Lawson Moving Lemma)
Suppose that $X\subset \P^n$ is a real projective variety of
dimension $m$. Fix nonnegative integers $r, s, e$ where $r+s\geq m$.
There exists a connected interval $I\subset \R$ containing 0 and a
map
$$\Psi: Z_r(X)\times I \longrightarrow Z_r(X)$$ such that
\begin{enumerate}
    \item $\Psi_0=$ identity.
    \item For $t\in I$, $\Psi_t$ is a real topological group homomorphism.
    \item For all cycles $Y, Z$ of dimensions $r, s$ and degree
    $\leq e$ and all $t\neq 0$, any component of excess
    dimension(i.e., $>r+s-m$) of $Y\bullet \Psi_t(Z)$ lies in the
    singular  locus of $X$ for any $t\neq 0$.
\end{enumerate}
\end{theorem}

\begin{proof}
We give a brief explanation why their proof of the moving lemma can
be reduced to real case. Let $F=(f_0, \cdots, f_m)$, a (m+1)-tuples
of real homogenous polynomials of degree $d$ in $\P^n$ where the
zero locus of $F$ is disjoint from $X$. Via the Veronese embedding,
we embed $\P^n$ into $\P^K$ where $K=\binom{n+d}{d}-1$. Let
$\pi_F:\P^M\cdots \rightarrow \P^m$ be the linear projection
determined by $F$ and $L(F)$ the center of $\pi_F$. $F$ determines a
finite map $p_F: X \longrightarrow \P^m$, see for example
\cite{Shafa}. For a cycle $Z$, $\pi^*_F(Z)=Z\#L(F)$,
$C_F(Z)=\pi^*_F(p_{F*}(Z))$, $R_F(Z)=C_F(Z)\bullet X-Z$ where $\#$
is the join. Since $F$ is real, $\pi^*_F$, $C_F$ and $R_F$ are real
maps. In \cite[Proposition 2.3]{FL2}, a map
$$\varphi_{L, D}(Z)=\pi_{0*}((L\# Z)\bullet D_1 \bullet \cdots
\bullet D_t)$$ is defined. We may take the effective divisors
$D_1,..., D_t$ to be real effective divisors, then $\varphi$ is a
real map, furthermore, the proof and the conclusion remain the same.

Friedlander and Lawson construct the map $\Psi$ from the map
$$\Psi_{\underline{N}}(Z,
p)=(-1)^{m+1}(M+1)\cdot
R_{F^*}(Z)+\sum_{i=1}^m(-1)^i\pi^*_{F^i}\{\Theta_{\underline{N},
p}\{p_{F^i*}(R_{F^{i-1}}\circ \cdots \circ R_{F^0}(Z))\}\}\bullet
X$$ for some $\underline{N}$ large enough where $\underline{N}$ is a
tuple of some positive integers, $M$ is a positive integer. The map
$\Theta_{\underline{N}, p}$ is constructed by iterating the map
$\varphi_{L, D}$ several times, so if we take $D$ to be a tuple of
real effective divisors, $\Theta_{\underline{N}, p}$ is also a real
map. It follows that $\Psi$ can be constructed to be a real map.
\end{proof}

Now we are able to give a proof of a duality theorem between reduced
real Lawson homology and reduced real morphic cohomology for
nonsingular real projective varieties.

\label{duality theorem}
\begin{theorem}
Suppose that $X, Y$ are nonsingular real projective varieties and
the dimension of $X$ is $m$. Then
\begin{enumerate}
    \item $R_r(Y)(X)$ is homotopy equivalent to $R_{r+m}(X\times Y)$
    \item $R^t(X)$ is homotopy equivalent to $R_m(X\times \A^t)$
    for any $t\geq 0$
    \item $R^t(X)$ is homotopy equivalent to $R_{m-t}(X)$ for $0\leq
    t\leq m$.
\end{enumerate}
Therefore, we have a group isomorphism
$$RL^tH^k(X)\cong RL_{m-t}H_{m-k}(X)$$ for $0\leq t, k \leq m$.
\end{theorem}

\begin{proof}
For each $e\geq 0$, let
$$K_e=\underset{d_1+d_2\leq e}
{\coprod}\mathscr{C}_{r, d_1}(X\times Y)\times \mathscr{C}_{r,
d_2}(X\times Y)/\sim$$ and let
$$K'_e=\underset{d_1+d_2\leq
e}{\coprod}\mathscr{C}_{r, d_1}(Y)(X)\times \mathscr{C}_{r,
d_2}(Y)(X)/\sim$$ where $\sim$ is the naive group completion
relation and $\mathscr{C}_{r, d}(Y)(X)=\mathscr{C}_{r+m, d}(X\times
Y)\cap \mathscr{C}_r(Y)(X)$.

Let $q$ be the quotient map from $Z_{r+m}(X\times Y)_{\R}$ to
$R_{r+m}(X\times Y)$ and let $q'$ be the quotient map from
$Z_r(Y)(X)_{\R}$ to $R_r(Y)(X)$. Let $\widetilde{K_e}=q(K_e\cap
Z_{r+m}(X\times Y)_{\R}), \widetilde{K'_e}=q'(K'_e\cap
Z_r(Y)(X)_{\R})$.

By the real Friedlander-Lawson moving lemma, we get two real maps
$$\Psi: Z_{r+m}(X\times Y)\times I \longrightarrow
Z_{r+m}(X\times Y),$$
$$\Psi': Z_r(Y)(X) \times I \longrightarrow Z_r(Y)(X)$$ where $I$
is some connected interval containing 0. By abuse of notation, we
will use the same notation to denote the two maps induced by $\Psi,
\Psi'$ respectively
$$\Psi: R_{r+m}(X\times Y)\times I \longrightarrow
R_{r+m}(X\times Y),$$
$$\Psi': R_r(Y)(X) \times I \longrightarrow R_r(Y)(X)$$

Restricting $\Psi, \Psi'$ to $\widetilde{K_e}\times I$ and
$\widetilde{K'_e}\times I$, we get two maps
$$R\phi_e=\widetilde{\Psi}|_{\widetilde{K}_e\times I}, \ \
R\phi'_e=\widetilde{\Psi'}|_{\widetilde{K'}_e\times I}.$$

The inclusion map $i:\mathscr{C}_r(Y)(X) \rightarrow
\mathscr{C}_{r+m}(X\times Y)$ induces the Friedlander-Lawson duality
map
$$\mathscr{D}:Z_r(Y)(X) \rightarrow Z_{r+m}(X\times Y)$$
As shown in \cite{FL3}, this map is continuous, injective but is not
a topological embedding.

Since $X, Y$ are real projective varieties, the duality map is real,
hence it induces an injective map
$$R\mathscr{D}: R_r(Y)(X) \longrightarrow R_{r+m}(X\times Y).$$
By Lemma \ref{compact subset}, the filtrations
$$\widetilde{K_0}\subset \widetilde{K_1} \subset \cdots =R_{r+m}(X\times
Y)$$ and
$$\widetilde{K'_0}\subset \widetilde{K'_1} \subset \cdots=R_r(Y)(X)$$
are locally compact and $R\mathscr{D}$ is filtration-preserving.

We have the following commutative diagrams:
$$\xymatrix{\widetilde{K'}_e \times I \ar[d]_{R\mathscr{D}\times Id} \ar[rr]^{R\phi'_e}&&
R_r(Y)(X) \ar[d]^{R\mathscr{D}}\\
\widetilde{K}_e \times I \ar[rr]^{R\phi_e} && R_{r+m}(X\times Y) }$$

$$\xymatrix{\widetilde{K'}_e \times\{1\} \ar[rr]^{R\phi'_e}
\ar[d]&&
R_r(Y)(X) \ar[d]\\
\widetilde{K}_e \times \{1\} \ar[rru]^{R\phi_e} \ar[rr]^{R\phi_e} &&
R_{r+m}(X \times Y)}
$$
and there is a map $\widetilde{\lambda}_e=R\phi_e$ from
$\widetilde{K_e}\times 1$ to $R_r(Y)(X)$. Thus by Lemma [5.2] in
\cite{FL3}, $R_r(Y)(X)$ is homotopy equivalent to $R_{r+m}(X\times
Y)$.

Furthermore, we have a commutative diagram of topological groups:
$$
\begin{array}{ccccc}
  R_0(\P^{t-1})(X) & \longrightarrow & R_0(\P^t)(X) & \longrightarrow & R^t(X) \\
  \downarrow & \ & \downarrow & \ & \downarrow \\
  R_m(X\times \P^{t-1}) & \longrightarrow & R_m(X\times \P^t) & \longrightarrow &
  R_m(X\times \A^t)\\
\end{array}
$$

The first two vertical arrows are homotopy equivalences which
implies the last one is also a homotopy equivalence. If $0\leq t
\leq m$, by the homotopy property of trivial bundle projection,
$R^t(X)$ is homotopy equivalent to $R_{m-t}(X)$.
\end{proof}

\begin{proposition}
Suppose that $X$ is a nonsingular real projective variety. Then
$R_0(X\times \A^t)$ is homotopy equivalent to $\Omega^{-t} R_0(X)$
where \index{$\Omega^{-t} R_0(X)$} $\Omega^{-t} R_0(X)$ is the
$t$-fold delooping of $R_0(X)$ given the infinite loop space
structure induced by the structure as a topological abelian group of
$R_0(X)$.
\end{proposition}

\begin{proof}
For two Eilenberg-Mac Lane spaces $K(G, i), K(H, j)$, denote
$$K(G, i)\bigotimes K(H, j)=K(G\otimes H, i+j)$$ and
$$(\prod^n_{i=1}(K(G_i, i)))\bigotimes (\prod^m_{j=1}(K(H_j,
j)))=\prod^{n+m}_{r=1}\prod_{i+j=r} K(G_i\otimes H_j, i+j).$$ From
Theorem A.5. in \cite{LLM1}, there is a canonical homotopy
equivalence between $R_0(X)$ and $\prod_{k\geq 0}K(H_k(ReX; \Z_2),
k)$. We will always consider the homotopic splitting of $R_0(X)$
into Eilenberg-Mac Lane spaces by this canonical homotopy
equivalence. By the K\"{u}nneth formula for $\Z_2$-coefficients and
the Dold-Thom theorem, we have $R_0(X\times
\P^{t-1})=R_0(X)\bigotimes R_0(\P^{t-1})$ and $R_0(X\times
\P^t)=R_0(X)\bigotimes R_0(\P^t)$. Since the inclusion map $\P^{t-1}
\hookrightarrow \P^t$ induces an isomorphism $i_*: H_k(\RP^{t-1};
\Z_2) \longrightarrow H_k(\RP^t; \Z_2)$ for $0\leq k < t$ in
homology, $i_*: \pi_k(R_0(\P^{t-1})) \longrightarrow
\pi_k(R_0(\P^t))$ is an isomorphism for $0\leq k<t$ by the Dold-Thom
theorem. Therefore $i_*: \pi_lR_0(X\times \P^{t-1})\longrightarrow
\pi_l R_0(X\times \P^t)$ is an isomorphism for $0\leq l <t$ and an
injection for $l\geq t$. From the long exact sequence
$$\cdots \longrightarrow \pi_lR_0(X \times \P^{t-1}) \longrightarrow
\pi_lR_0(X \times \P^t) \longrightarrow \pi_l(R_0(X\times
\A^t))\longrightarrow \cdots$$ we see that $\pi_l(R_0(X\times
\A^t))=0 \mbox{ if } 0\leq l<t$ and $\pi_l(R_0(X\times
\A^t))=\frac{\pi_lR_0(X \times \P^t)}{\pi_lR_0(X \times \P^{t-1})}$
if $l\geq t$. Use the canonical homotopical splitting of $R_0(X)$
into a product of Eilenberg-Mac Lane spaces and denote the $i$-th
component of the Eilenberg-Mac Lane space of $R_0(X)$ to be $R_{0,
i}(X)$. Then for $l\geq t$, from above calculation, we have $R_{0,
l}(X\times \A^t)=R_{0, l-t}(X)\bigotimes R_{0, t}(\P^t)=R_{0,
l-t}(X)\bigotimes K(\Z_2, t)=\Omega^{-t} R_{0, l-t}(X)$. Thus
$R_0(X\times \A^t)=\Omega^{ -t}R_0(X)$.
\end{proof}

By taking homotopy groups, we have the following result:
\begin{corollary}
For a nonsingular real projective variety $X$,
$$RL_0H_k(X\times\A^t)=\left\{%
\begin{array}{ll}
    RL_0H_{k-t}(X), & \hbox{ if } k\geq t; \\
    0, & \hbox{ if } k<t. \\
\end{array}%
\right.$$
\end{corollary}

\begin{corollary}
Suppose that $X$ is a nonsingular real projective variety of
dimension $m$ and its real points $ReX$ is connected and of real
dimension $m$. Then
$$RL^tH^k(X)=RL^mH^k(X)=H^k(ReX; \Z_2)$$ for
$t\geq m \geq k$ and $RL^tH^k(X)=0$ for $t\geq k>m$.
\end{corollary}

\begin{proof}
$RL^tH^k(X)=\pi_{t-k}R^t(X)=\pi_{t-k}R_m(X\times
\A^t)=\pi_{t-k}R_0(X\times \A^{t-m})=\pi_{t-k}\Omega^{-(t-m)}R_0(X)$ $$=\left\{%
\begin{array}{ll}
    0, & \hbox{ if } k>m \\
    \pi_{m-k}R_0(X)=H_{m-k}(ReX; \Z_2)=H^k(ReX; \Z_2)=RL^mH^k(X) & \hbox{ if } k\leq m \\
\end{array}%
\right.$$
\end{proof}

\begin{corollary}
$R^t(\P^n)=$
$\left\{%
\begin{array}{ll}
    \prod^{t}_{i=0}K(\Z_2, i) & \hbox{ if } t\leq n \\
    \prod^n_{i=0}K(\Z_2, i+t-n) & \hbox{ if } t>n \\
\end{array}
\right. $
\end{corollary}

\begin{proof}
If $t\leq n$, by the Duality Theorem, $R^t(\P^n)=R_n(\P^n\times
\A^t)=R_{n-t}(\P^n)=\prod^{t}_{i=0}K(\Z_2, i)$; if $t>n$, then
$R^t(\P^n)=R_n(\P^n\times \A^t)=R_0(\P^n\times
\A^{t-n})=\Omega^{-(t-n)}R_0(\P^n)=\Omega^{-(t-n)}
\prod^n_{i=0}K(\Z_2, i)=\prod^n_{i=0}K(\Z_2, i+t-n)$.
\end{proof}

\section{Operations and Maps}

\subsection{Natural Transformations From Reduced Real Morphic Cohomology To Singular Cohomology}
Suppose that $X, Y$ are real projective varieties. It is a general
fact that for a topological group $G$ and a locally compact
Hausdorff space $A$, $Map(A, G)$, the space of all continuous
functions from $A$ to $G$ with the compact-open topology, is a
topological group under point-wise multiplication. By using the
graphing construction of Friedlander and Lawson, there is a natural
transformation $j:Z_r(Y)(X) \rightarrow Map(X, Z_r(Y))$ induced by
the inclusion map $i:\mathscr{M}or(X, \mathscr{C}_r(Y)) \rightarrow
Map(X, Z_r(Y))$. The restriction of $j$ to real cocycles gives us a
map $j:Z_r(Y)(X)_{\R} \rightarrow Map(X, Z_r(Y))$. For $f\in
Z_r(Y)(X)_{\R}$, if we restrict $f$ to $Re(X)$, the set of real
points of $X$, since $f(x)=\overline{f(\bar{x})}=\overline{f(x)}$,
the image of $f$ lies in $Z_r(Y)_{\R}$. Composing the map $j$ with
the restriction map, we have a continuous map $$\Psi':
Z_r(Y)(X)_{\R} \longrightarrow Map(ReX, Z_r(Y)_{\R}).$$ Composing
again with the quotient map $q: Z_r(Y)_{\R} \longrightarrow R_r(Y)$,
we have a continuous map
$$\Psi'':Z_r(Y)(X)_{\R}\longrightarrow Map(ReX, R_r(Y)).$$
If $f=g+\overline{g}$,
$\Psi''(f)(x)=g(x)+\overline{g(\bar{x})}+Z_r(Y)^{av}=g(x)+\overline{g(x)}+
Z_r(Y)^{av}=Z_r(Y)^{av}$, so $Z_r(Y)(X)^{av} \subset Ker\Psi''$,
therefore $\Psi''$ induces a continuous map
$$\Psi: R_r(Y)(X) \longrightarrow Map(ReX, R_r(Y)).$$

We summarize the above construction in the following diagram:
$$\xymatrix{ Z_r(Y)(X)_{\R} \ar[dd] \ar[rrr]^j \ar[rrrd]^{\Psi'}
\ar[rrrdd]^{\Psi''}
& & & Map(X, Z_r(Y)) \ar[d]^{\mbox{ restriction }}\\
& & & Map(ReX, Z_r(Y)_{\R}) \ar[d]^q\\
R_r(Y)(X) \ar[rrr]^{\Psi} & & & Map(ReX, R_r(Y)) }$$

Now we are going to construct a map from $R^t(X)$ to $Map(ReX,
R_0(\A^t))$.

Let $q': R_0(\P^t) \longrightarrow R_0(\A^t)$ be the quotient map.
Then $q'$ induces a map, also denoted as $q'$, from $Map(ReX,
R_0(\P^t))$ to $Map(ReX, R_0(\A^t))$. Let $\Phi=q'\circ \Psi:
R_0(\P^t)(X) \longrightarrow Map(ReX, R_0(\A^t))$. For $f\in
R_0(\P^{t-1})(X)$, $f\in ker(\Phi)$, so $\Phi$ induces a map
$\Phi^t:R^t(X) \longrightarrow Map(ReX, R_0(\A^t))$. We summarize
the construction of $\Phi^t$ in following:
$$\Phi^t(f+Z_0(\P^t)(X)^{av}+R_0(\P^{t-1})(X))=q'\circ q\circ (f|_{Re
X}).$$

Since $R_0(\A^t)$ is the Eilenberg-Mac Lane space $K(\Z_2, t)$,
taking homotopy groups, $\Phi^t$ induces a group homomorphism:
$$\Phi^{t, k}: RL^tH^k(X) \longrightarrow H^k(ReX; \Z_2).$$

\begin{proposition}
$\Phi^{t,k}: RL^tH^k(X) \longrightarrow H^k(ReX; \Z_2)$ is a natural
transformation for each $k$ with $0\leq k \leq t$.
\end{proposition}

\begin{proof}
Suppose that $f:X \longrightarrow Y$ is a morphism between two real
projective varieties. We have a pullback $f^*:R^t(Y) \longrightarrow
R^t(X)$ (see Proposition \ref{pullback}) by mapping
$\phi+R_0(\P^{t-1})(Y)$ to $\phi \circ f+ R_0(\P^{t-1})(X)$. Since
$f$ is a real map, it maps $Re(X)$ to $Re(Y)$, thus $f$ induces a
map $f^*: Map(ReY, R_0(\A^t)) \longrightarrow Map(ReX, R_0(\A^t))$.
It is easy to check that the following diagram commutes:

$$\xymatrix{
R^t(Y) \ar[r]^{f^*} \ar[d]_{\Phi^t} & R^t(X) \ar[d]^{\Phi^t} \\
Map(ReY, R_0(\A^t)) \ar[r]^{f^*} & Map(ReX, R_0(\A^t))}$$ which then
induces a commutative diagram of homotopy groups. For $g:Y
\rightarrow Z$ another morphism between real projective varieties,
we have $\Phi^{t, k}(g\circ f)=\Phi^{t, k}(f)\circ \Phi^{t, k}(g)$.
\end{proof}

\subsection{Cup Product}
Suppose that $X, Y$ are two real projective varieties. The join
$X\#Y$ of $X, Y$ is again real. If $[x_0:...:x_n:y_0:...:y_m]\in X\#
Y$, the conjugation on $X\#Y$ is defined by
$\overline{[x_0:...:x_n:y_0:...:y_m]}=[\bar{x}_0:...:\bar{x}_n:\bar{y}_0:...:\bar{y}_m]$.

From Section 6 in \cite{FL1}, there is a continuous pairing induced
by the join of varieties:
$$Z^s(Y)(X)\wedge Z^{s'}(Y')(X') \overset{\#}{\longrightarrow}
Z^{s+s'}(Y\#Y')(X\times X')$$ given by
$$(\varphi\#\varphi')(x, x')=(\varphi(x))\#(\varphi'(x'))$$

It is easy to see that this pairing is real. Therefore, for real
projective varieties $X, X', Y, Y'$, $\#$ reduces to a continuous
pairing
$$R^s(Y)(X)\wedge R^{s'}(Y')(X') \overset{\#}{\longrightarrow}
R^{s+s'}(Y\#Y')(X\times X').$$

Taking homotopy groups, we have a pairing:
$$RL^sH^q(X; Y)\otimes RL^{s'}H^{q'}(X'; Y') \longrightarrow
RL^{s+s'}H^{q+q'}(X\times X'; Y\#Y'),$$ and when restricted to the
diagonal of $X\times X$, it determines a cup product:
$$RL^sH^q(X; Y)\otimes RL^{s'}H^{q'}(X; Y')
\overset{\#}{\longrightarrow} RL^{s+s'}H^{q+q'}(X; Y\#Y').$$

Now take $Y=\P^s, Y'=\P^{s'}$, the get a pairing
$$R_0(\P^s)(X)\wedge R_0(\P^{s'})(X) \longrightarrow
R_1(\P^{s+s'+1})(X\times X)$$ Since
$$R_0(\P^{s-1})(X)\wedge R_0(\P^{s'})(X)\rightarrow
R_1(\P^{s+s'})(X)$$ and
$$R_0(\P^s)(X)\wedge R_0(\P^{s'-1})(X) \rightarrow
R_1(\P^{s+s'})(X),$$ the pairing reduces to a pairing on the reduced
real cocycle groups :
$$R^s(X)\wedge R^{s'}(X) \longrightarrow \frac{R_1(\P^{s+s'+1})(X)}{R_1(\P^{s+s'})(X)}
\overset{\susp^{-1}}{\longrightarrow} R^{s+s'}(X\times X)$$ By
restricting to the diagonal of $X\times X$, we get a pairing
$$R^s(X)\wedge R^{s'}(X) \longrightarrow R^{s+s'}(X)$$ which gives
us a commutative cup product in reduced real morphic cohomology:
$$RL^sH^q(X)\otimes RL^{s'}H^{q'}(X) \longrightarrow
RL^{s+s'}H^{q+q'}(X).$$

\begin{proposition}
Suppose that $X, X', Y, Y', W, W', Z, Z'$ are all real projective
varieties and $f: X \longrightarrow X', g: W \longrightarrow W'$ are
real morphisms. Then we have the following commutative diagrams:
\begin{enumerate}
\item $$\xymatrix{R^s(Y)(X')\wedge R^{s'}(Z)(W') \ar[r]^-{\#}
\ar[d]_{f^*\wedge g^*} & R^{s+s'}(Y\# Z)(X'\times W') \ar[d]^{(f\times g)^*}\\
R^s(Y)(X)\wedge R^{s'}(Z)(W) \ar[r]^-{\#} & R^{s+s'}(Y\# Z)(X\times
W)}$$

\item $$\xymatrix{R^s(X')\wedge R^{s'}(X') \ar[r]^{\ \ \ \#}
\ar[d]_{f^*\wedge f^*} & R^{s+s'}(X')
\ar[d]^{f^*}\\
R^s(X)\wedge R^{s'}(X) \ar[r]^{\ \ \ \#} & R^{s+s'}(X)}$$

\item $$\xymatrix{R^s(X)(Y)\wedge R^{s'}(W)(Z) \ar[r]^-{\#}
\ar[d]_{f_*\wedge g_*} & R^{s+s'}(X\# W)(Y\times Z) \ar[d]^{(f\times g)_*}\\
R^{s-c}(X')(Y)\wedge R^{s'-c'}(W')(Z) \ar[r]^-{\#} &
R^{(s+s')-(c+c')}(X'\# W')(Y\times Z)}$$ where $c=dimX-dimX'$,
$c'=dimW-dimW'$.

In short, cup product in reduced real morphic cohomology is natural
with respect to real morphisms.
\end{enumerate}
\end{proposition}

\begin{proof}
We observe that $$(f^*\varphi)\#(g^*\varphi')=(\varphi\circ f)\#
(\varphi' \circ g)=(\varphi \# \varphi')\circ (f\times g)=(f\times
g)^*(\varphi \# \varphi'),$$
$$(f_*\alpha)\#(g_*\alpha')=(f\times g)_*(\alpha\# \alpha')$$ and
then check that they pass to reduced real cocycles.
\end{proof}

From Lawson suspension theorem, we obtain a canonical homotopy
equivalence:
$$Z^q(\P^n)\longrightarrow \prod^q \overset{\mbox{def}}{=}
\prod^q_{k=0}K(\Z, 2k)$$ for all $n\geq q$. In \cite{LM}, Lawson and
Michelsohn showed that the complex join
$$\#_{\C}:\overset{q}{\prod}\times \overset{q'}{\prod}
\longrightarrow \overset{q+q'}{\prod}$$ has the property that
$$\#_{\C}^*(\iota_{2k})=\sum_{\substack{r+s=k}}\iota_{2r}\otimes \iota_{2s}$$
in integral cohomology  where $\iota_{2k}$ is the generator of
$$H^{2k}(K(\Z, 2k); \Z)\cong \Z.$$
The cup product and cross product in cohomology are characterized by
some axioms which can be found for example in \cite{AGP}. From the
above result, it is not difficult to check that the complex join
induces the cross product in integral cohomology. Following a
similar approach, Lam in \cite{Lam} showed that the corresponding
result holds for the pairing:
$$\#_{\R}: \overset{q}{\prod}_{\R}\times \overset{q'}{\prod}_{\R}
\longrightarrow \overset{q+q'}{\prod}_{\R}$$ where
$\overset{q}{\prod}_{\R}=K(\Z_2, 0)\times K(\Z_2, 1)\times \cdots
\times K(\Z_2, q)$, thus $\#$ induces the cross product in
$\Z_2$-cohomology.

We now show that the natural transformation from reduced real
morphic cohomology to $\Z_2$-singular cohomology is a ring
homomorphism.

\begin{proposition}
Suppose that $X$ and $X'$ are real projective varieties. Then for
all $s, s'$ and $q, q'$, there are commutative diagrams:
\begin{enumerate}
\item
$$\xymatrix{RL^sH^q(X)\otimes RL^{s'}H^{q'}(X') \ar[d]_{\Phi^{s, q}\otimes
\Phi^{s', q'}} \ar[r]^{\#} &
RL^{s+s'}H^{q+q'}(X\times X') \ar[d]^{\Phi^{s+s', q+q'}}\\
H^q(ReX; \Z_2)\otimes H^{q'}(ReX'; \Z_2) \ar[r]^-{\times} &
H^{q+q'}(Re(X\times X'); \Z_2)}$$ where the lower horizontal arrow
is the usual cross product in $\Z_2$-coefficients.

\item
$$\xymatrix{RL^sH^q(X)\otimes RL^{s'}H^{q'}(X) \ar[d]_{\Phi^{s, q}\otimes
\Phi^{s', q'}} \ar[r]^-{\#} &
RL^{s+s'}H^{q+q'}(X) \ar[d]^{\Phi^{s+s', q+q'}}\\
H^q(ReX; \Z_2)\otimes H^{q'}(ReX; \Z_2) \ar[r]^-{\smile} &
H^{q+q'}(ReX; \Z_2)}$$ where $\smile$ is the cup product in singular
cohomology with $\Z_2$-coefficients.
\end{enumerate}

In particular, the map
$$\Phi^{*, *}: RL^*H^*(X)\longrightarrow H^*(ReX; \Z_2)$$
is a graded-ring homomorphism.
\end{proposition}

\begin{proof}
Consider the following commutative diagram:
$$\xymatrix{R^s(X)\wedge R^{s'}(X') \ar[r]^-{\#} \ar[d] &
R^{s+s'}(X\times X') \ar[d]\\
Map(ReX, R_0(\A^s))\wedge Map(ReX', R_0(\A^{s'})) \ar[r]^-{\#} &
Map(Re(X\times X'), R_0(\A^{s+s'}))}$$ Since the lower arrow in the
diagram above is exactly the map classifying the cross product in
$\Z_2$-cohomology, taking homotopy groups, proves the first
assertion. Let $X'=X$ and $\Delta:X\longrightarrow X\times X$ be the
diagonal map. Since $x\smile y=\Delta^*(x\times y)$ for $x, y \in
H^q(ReX; \Z_2)$, we complete the proof by composing the first
diagram with the following commutative diagram:
$$\xymatrix{RL^{s+s'}H^{q+q'}(X\times X) \ar[r]^-{\Delta^*} \ar[d]_{\Phi^{s+s', q+q'}}&
RL^{s+s'}H^{q+q'}(X) \ar[d]^{\Phi^{s+s', q+q'}}\\
H^{q+q'}(ReX\times ReX; \Z_2) \ar[r]^-{\Delta^*} & H^{q+q'}(ReX;
\Z_2)}$$
\end{proof}

\subsection{The $S$-Map In Reduced Real Morphic Cohomology}
Let us give a construction of the $S$-map in reduced real morphic
cohomology parallel to the construction of Lima-Filho in Lawson
homology. Fix $x_{\infty} \in \RP^1$ and consider the following
diagram
$$\xymatrix{R^t(X)\times \RP^1 \ar[rr] \ar[d] & & R^t(X) \times R_0(\P^1)
\ar[d]^{\#}\\
R^t(X)\wedge \RP^1  \ar[rr]^{\eta} \ar[rrd]^{S} & &
\frac{R_1(\P^{t+2})(X)}{R_1(\P^{t+1})(X)}
\ar[d]^{\susp^{-1}}\\
& & \frac{R_0(\P^{t+1})(X)}{R_0(\P^t)(X)} \ar[d]^{\parallel}\\
& & R^{t+1}(X)}$$

The top horizontal row sends $(f+R_0(\P^{t-1})(X), x)\in
R^t(X)\times \RP^1$ to $(f+R_0(\P^{t-1})(X),
x-x_{\infty}+Z_0(\P^1)^{av})$ and the join on the right hand side
sends $(f+R_0(\P^{t-1})(X), y)$ to $f\#y+R_1(\P^{t+1})(X)$ where
$f\#y\in R_1(\P^{t+2})(X)$. This map reduces to map $\eta$ from
$R^t(X)\wedge \RP^1$ to $\frac{R_1(\P^{t+2})(X)}{R_1(\P^{t+1})(X)}$.
Define the $S$-map to be $\susp^{-1} \circ \eta$ where $\susp^{-1}$
is the homotopy inverse of the suspension map.

So we have a map from $R^t(X)\wedge \S^1 \longrightarrow R^{t+1}(X)$
which induces a map in the reduced real morphic cohomology:
$$S: RL^tH^k(X) \longrightarrow RL^{t+1}H^k(X).$$

\begin{proposition}
For a real projective variety $X$, the following diagram commutes:
$$\xymatrix{RL^tH^k(X) \ar[rr]^S \ar[rd]_{\Phi^{t, k}}& &
\ar[ld]^{\Phi^{t+1, k}} RL^{t+1}H^k(X)\\
& H^k(ReX; \Z_2) &}$$ \label{s map commutes}
\end{proposition}

\begin{proof}
Fix a point $x_{\infty}\in \RP^1$. Consider the following diagram:
$$\xymatrix{R^t(X)\wedge \RP^1 \ar[r]^{\eta} \ar[d]_{\Phi^t\times Id} &
\frac{R_1(\P^{t+2})(X)}{R_1(\P^{t+1})(X)} \ar[d] & \ar[l]_{\susp}
R^{t+1}(X)
\ar[d]^{\Phi^{t+1}}\\
Map(ReX, R_0(\A^t))\wedge \RP^1 \ar[r] & Map(ReX,
\frac{R_1(\P^{t+2})}{R_1(\P^{t+1})}) & \ar[l]_{\susp} Map(ReX,
R_0(\A^{t+1}))}$$ The right arrow in the bottom row sends $ f \wedge
x$ to the map $(f\wedge
x)(y)=f(y)\#((x-x_{\infty})+Z_0(\P^1)^{av})$. It is easy to see that
these two squares commute and then take the homotopy groups.
\end{proof}

\begin{corollary}
For a real projective variety $X$, define $\mathscr{F}^t=\Phi^{t,
k}(RL^tH^k(X))$. By the above Proposition, we have a filtration:
$$\mathscr{F}^k \subseteq \mathscr{F}^{k+1} \subseteq \cdots
\subseteq H^k(ReX; \Z_2)$$ which is analogous to the
\textquotedblleft topological filtration" defined by Friedlander and
Mazur in \cite{FM}.
\end{corollary}

\subsection{The H-operator In Reduced Real Morphic Cohomology}
From the construction of $S$-map, we have a map
$$H: R^t(X)\hookrightarrow R^t(X)\times R_0(\P^1) \longrightarrow
R^{t+1}(X)$$ and then taking the homotopy groups, we have
$$H: RL^tH^k(X) \longrightarrow RL^{t+1}H^{k+1}(X)$$
which is the $H$-operator in reduced real morphic cohomology.

\subsection{Natural Transformations From Reduced Real Lawson Homology To Singular Homology}
In Lawson homology, there is a natural transformation
$$\Phi_{r, k}: L_rH_k(X) \longrightarrow H_k(X; \Z)$$ defined by
iterating the $s$-map in Lawson homology $r$ times and then
composite with the Dold-Thom isomorphism. We show that similar
construction is valid in reduced real Lawson homology and we get a
natural transformation from reduced real Lawson homology to singular
homology of the real points with $\Z_2$-coefficients.

Fix a point $x_{\infty}\in \RP^1$.
$$\xymatrix{R_t(X)\times \RP^1 \ar[r] \ar[d] & R_t(X) \times R_0(\P^1)
\ar[d]^{\#}\\
R_t(X)\wedge \RP^1  \ar[r]^{\eta} \ar[rd]^{s} & R_{t+1}(X\# \P^1) \ar[d]^{\susp^{-2}}\\
& R_{t-1}(X)}$$ The map on the top horizontal row sends
$(V+Z_t(X)^{av}, x)\in R_t(X)\times \RP^1$ to $(V+Z_t(X)^{av},
x-x_{\infty}+Z_0(\P^1)^{av})\in R_t(X)\times R_0(\P^1)$. It reduces
to a map $\eta$ from the smash product of $R_t(X)$ and $\RP^1$ to
$R_{t+1}(X\# \P^1)$. Finally, we take the homotopy inverse of the
suspension map twice. So we get a map $s: \S^1\wedge R_t(X)
\longrightarrow R_{t-1}(X)$ which induces a map in the reduced real
Lawson homology
$$s: RL_tH_k(X)\longrightarrow RL_{t-1}H_{k}(X).$$
We iterate this map $t$ times and then apply the Dold-Thom
isomorphism $\tau:\pi_kR_0(X) \overset{\cong}{\longrightarrow}
H_k(ReX; \Z_2)$. This gives us a map
$$\Phi_{t, k}: RL_tH_k(X)\longrightarrow H_k(ReX; \Z_2)$$ where
$\Phi_{t, k}=\tau\circ s^t$.

\begin{proposition}
For a morphism $f:X\longrightarrow Y$ between real projective
varieties $X, Y$, the following diagram commutes:
$$\xymatrix{RL_rH_k(X) \ar[r]^{\Phi_{r, k}} \ar[d]_{f_*} & H_k(ReX; \Z_2) \ar[d]^{f_*}\\
RL_rH_k(Y) \ar[r]^{\Phi_{r, k}} & H_k(ReY; \Z_2)}$$ Thus $\Phi$ is a
natural transformation from reduced real Lawson homology to singular
homology of real points with $\Z_2$-coefficients.
\end{proposition}

\begin{proof}
The following diagram commutes:
$$\xymatrix{R_r(X)\wedge \RP^1 \ar[r]^{\eta} \ar[d]^{f_*\wedge Id} & R_{r+1}(X\#\P^1)
 \ar[d]^{(f\#Id)_*} & R_{r-1}(X) \ar[l]_-{\susp ^2} \ar[d]^{f_*}\\
R_r(Y)\wedge \RP^1 \ar[r]^{\eta} & R_{r+1}(Y\#\P^1) &
\ar[l]_-{\susp^2} R_{r-1}(Y)}$$

Taking homotopy groups, we have the following commutative diagram of
$s$-map
$$\xymatrix{RL_rH_k(X) \ar[r]^s \ar[d]^{f_*} & RL_{r-1}H_k(X) \ar[d]^{f_*} \\
RL_rH_k(Y) \ar[r]^s & RL_{r-1}H_k(Y)}$$

And by the functoriality of Dold-Thom isomorphism, we have the
commutative diagram:
$$\xymatrix{RL_0H_k(X) \ar[r]^{\tau} \ar[d]^{f_*} & H_k(ReX; \Z_2) \ar[d]^{f_*}\\
RL_0H_k(Y) \ar[r]^{\tau} & H_k(ReY; \Z_2)}$$

Applying the commutative diagram of $s$-map $r$ times and applying
the commutative diagram of Dold-Thom isomorphism, we get the
required commutativity.
\end{proof}

\subsection{Filtrations}
\begin{definition}(The geometric filtration)
Let $X$ be a real projective variety and denote by
$RG_jH_n(X)\subset H_n(ReX; \Z_2)$ the subspace of $H_n(ReX; \Z_2)$
generated by the image of maps $H_n(ReY; \Z_2)\longrightarrow
H_n(ReX; \Z_2)$ induced from all morphisms $Y\longrightarrow X$ of
real projective variety $Y$ of dimension $\leq 2n-j$. The subspaces
$RG_jH_n(X)$ form a decreasing filtration:
$$\cdots \subset RG_jH_n(X) \subset RG_{j-1}H_n(X) \subset \cdots \subset
RG_0H_n(X)\subset H_n(ReX; \Z_2)$$ which is called the geometric
filtration.
\end{definition}

The $s$-map in reduced real Lawson homology enables us to define a
filtration which is analogous to the topological filtration in
Lawson homology.

\begin{definition}(The topological filtration)
Suppose that $X$ is a real projective variety. Let $RT_tH_n(X)$
denote the subspace of $H_n(ReX; \Z_2)$ given by the image of
$\Phi_{t, n}$, i.e.,
$$RT_tH_n(X)=\Phi_{t, n}(RL_tH_n(X)).$$
The subspaces $RT_tH_n(X)$ form a decreasing filtration:
$$\cdots \subset RT_tH_n(X) \subset RT_{t-1}H_n(X) \subset \cdots
\subset RT_0H_n(X)=H_n(ReX; \Z_2),$$ and $RT_tH_n(X)$ vanishes if
$t>n$. This filtration is called the topological filtration.
\end{definition}

It was conjectured by Friedlander and Mazur in \cite{FM} that the
topological filtration and geometric filtration in Lawson homology
coincide. We post the similar conjecture for reduced real Lawson
homology.

\begin{conjecture}
For a smooth real projective variety $X$, the topological filtration
and the geometric filtration in reduced real Lawson homology
coincide, i.e.,
$$RT_tH_n(X)=RG_tH_n(X).$$
\end{conjecture}

\subsection{The h-Operator }
Let us now construct a map which is the analogue of the $h$-operator
in Lawson homology. Consider the map
$$h: R_r(X)\hookrightarrow R_r(X) \times R_0(\P^1)
\overset{\#}{\longrightarrow} R_{r+1}(X\#\P^1)
\overset{\susp^{-2}}{\longrightarrow} R_{r-1}(X)$$ which induces a
map, called the $h$-operator, in reduced real Lawson homology:
$$h: RL_rH_k(X)\longrightarrow RL_{r-1}H_{k-1}(X).$$

\subsection{Slant Product} \label{slant product} Suppose that $X, Y$
are projective varieties. Define a product
$$Z_r(Y)(X)\times Z_p(X) \longrightarrow Z_{r+p}(Y)$$
by sending $(f, V)$ to $Pr_*(f\bullet (V\times Y))$ where
$Pr:X\times Y \longrightarrow Y$ is the projection and we consider
$f$ as a cycle in $X\times Y$ which intersects $V\times Y$ properly.
It is proved in \cite[Proposition 7.1]{FL1} that this is a
continuous pairing.

Now suppose that $X, Y$ are real projective varieties. Since in the
definition, each operation is real, the above pairing reduces to a
pairing
$$R_r(Y)(X)\times R_p(X)\longrightarrow R_{r+p}(Y)$$
and it is easy to see that it reduces again to
$$R_r(Y)(X)\wedge R_p(X)\longrightarrow R_{r+p}(Y).$$
Therefore, we have a slant product:
$$RL^rH^k(X; Y)\otimes RL_pH_n(X) \longrightarrow
RL_{r+p}H_{n-k+2r}(Y).$$

Fix $r=0, p\leq t$. For $Y=\P^{t-1}$, the product from the above
construction sends
$$R_0(\P^{t-1})(X)\wedge R_p(X) \longrightarrow R_p(\P^{t-1}).$$
Therefore, the product
$$R_0(\P^t)(X)\wedge R_p(X) \longrightarrow R_p(\P^t)$$
reduces to
$$\frac{R_0(\P^t)(X)}{R_0(\P^{t-1})(X)}\wedge
R_p(X)\longrightarrow \frac{R_p(\P^t)}{R_p(\P^{t-1})}.$$ Since
$\frac{R_p(\P^t)}{R_p(\P^{t-1})}=R_p(\A^t)\cong R_0(\A^{t-p})$, we
have a pairing:
$$R^t(X)\wedge R_p(X) \longrightarrow R_0(\A^{t-p})$$ which
induces a Kronecker pairing:
$$RL^tH^k(X)\otimes RL_pH_k(X)\longrightarrow \Z_2$$ for $p\leq
k\leq t$.

\begin{proposition}
Let $X$ be a real projective variety. Then for all $s \geq q$, the
diagram
$$\xymatrix{RL^sH^q(X)\otimes RL_0H_q(X) \ar[dd]_{\Phi^{s, q}\otimes \Phi_{0, q}}
 \ar[rd]^{\kappa} & \\
& \Z_2\\
H^q(ReX; \Z_2)\otimes H_q(ReX; \Z_2) \ar[ru]^{\kappa^{top}} &}$$
commutes; i.e., under the natural transformation to $\Z_2$-singular
theory, the Kronecker pairing introduced above is carried to the
topological one. \label{kronecker pairing}
\end{proposition}

\begin{proof}
The canonical homeomorphism between $R_0(X)$ and
$\frac{Z_0(ReX)}{2Z_0(ReX)}$ is given by sending an element $c\in
R_0(X)$ to the cycle formed by real points of $X$ that represent $c$
which is unique modulo 2, i.e., there exist unique $x_1, ...., x_k
\in ReX$ such that $c=\sum_{i=1}^kx_i+Z_0(X)^{av}$ and $c$ is mapped
to $\sum_{i=1}^kx_i+2Z_0(ReX)$ as the proof in Proposition \ref{real
points isomorphic}. In the following, we will assume that $c$ is
represented by real points.

Consider the following diagram:
$$\xymatrix{R^s(X)\wedge R_0(X) \ar[dd]_{\Phi^s\wedge Id} \ar[rd] & \\
& R_0(\A^s)\\
Map(ReX, R_0(\A^s))\wedge R_0(X) \ar[ru] &}$$

The slant product on the top horizontal arrow sends
$(f+Z_0(\P^s)^{av}+R_0(\P^{s-1})(X), \sum_i x_i+Z_0(X)^{av})$ to
$(\sum_if(x_i)+Z_0(\P^s)^{av}+R_0(\P^{s-1}))$ and the pairing in the
bottom arrow sends $(\varphi, \sum_ix_i)$ to $\sum_i\varphi(x_i)$.
Under the natural maps, it is not difficult to see that the diagram
commutes.

Taking homotopy groups, from the bottom horizontal arrow, we get a
pairing:
$$\kappa: \xymatrix{H^q(ReX; \Z_2)\otimes H_q(ReX; \Z_2) \ar[r] & \Z_2}$$
Observe that the following diagram commutes:
$$\xymatrix{Map(ReX, (\frac{Z_0(\S^t)}{2Z_0(\S^t)})^0)\wedge
\frac{Z_0(ReX)}{2Z_0(ReX)} \ar[d] \ar[r]
&(\frac{Z_0(\S^t)}{2Z_0(\S^t)})^0 \ar[d]\\
Map(ReX, R_0(\A^t))\wedge R_0(X)  \ar[r] & R_0(\A^t)}$$ where
$(\frac{Z_0(\S^t)}{2Z_0(\S^t)})^0=K(\Z_2, t)$ is the path-connected
component of $\frac{Z_0(\S^t)}{2Z_0(\S^t)}$ which contains the
identity and the canonical quotient map $S^t\longrightarrow \RP^t$
induces a map from $(\frac{Z_0(\S^t)}{2Z_0(\S^t)})^0$ to $R_0(\A^t)$
and hence a map from $Map(ReX, (\frac{Z_0(\S^t)}{2Z_0(\S^t)})^0)$ to
$Map(ReX, R_0(\A^t))$. Each vertical arrow of the diagram above is a
homotopy equivalence. Thus to prove this Proposition it will suffice
to establish the following lemma which is a generalization of Lemma
8.3 in \cite{FL1}.
\end{proof}

We use the convention that $\Z_0=\Z$ and
$\frac{Z_0(Y)}{0Z_0(Y)}=Z_0(Y).$
\begin{lemma}
For any finite CW-complex $Y$ and $p=0$ or a prime number, the
pairing
$$Map\left(Y, (\frac{Z_0(\S^t)}{pZ_0(S^t)})^0\right)\wedge
\frac{Z_0(Y)}{pZ_0(Y)} \longrightarrow
\left(\frac{Z_0(\S^t)}{pZ_0(S^t)}\right)^0$$ sending $(f,
\sum_in_iy_i+pZ_0(Y))$ to $\sum_in_if(y_i)+pZ_0(\S^t)$ induces a
pairing
$$\kappa: H^q(Y; \Z_p)\otimes H_q(Y; \Z_p) \longrightarrow \Z_p$$
which is the topological Kronecker pairing $\kappa^{top}$ where
$K=(\frac{Z_0(\S^t)}{pZ_0(\S^t)})^0=K(\Z_p, t)$ is the component of
$\frac{Z_0(\S^t)}{pZ_0(\S^t)}$ containing the identity.
\end{lemma}

\begin{proof}
For a continuous map $\varphi: Y \longrightarrow Z$, from the
definition, $\kappa$ has following naturality property:

\begin{equation} \kappa(\varphi^*\alpha, u)=\kappa(\alpha, \varphi_*u)
\label{naturality of kronecker pairing}
\end{equation}for all $\alpha\in H^q(Z; \Z_p)$ and $u\in H_q(Y; \Z_p)$.

Let $Y=\S^q$, $K=(\frac{Z_0(\S^t)}{pZ_0(\S^t)})^0=K(\Z_p, t)$. A
generator of $H^q(\S^q; \Z_p)=\pi_{t-q}Map(\S^q, K)$ is given by the
homotopy class of $f:\S^{t-q}\longrightarrow Map(\S^q, K)$ where
$f(x)(y)=x\wedge y+pZ_0(\S^t)$ and a generator of $H_q(\S^q;
\Z_p)=\pi_q(\frac{Z_0(\S^q)}{pZ_0(\S^q)})$ is given by the homotopy
class of $g: \S^q \longrightarrow \frac{Z_0(\S^q)}{pZ_0(\S^q)}$
defined by $g(x)=x+pZ_0(\S^q)$. Then the map $f\wedge g:
\S^{t-q}\wedge \S^q \longrightarrow K$ given by
$$(f\wedge g)(x\wedge y)=f(x)(g(y))=f(x)(y+pZ_0(S^t))=x\wedge y+ pZ_0(\S^t)$$
is the generator of $H_t(K; \Z_p)=\Z_p$. Thus $\kappa=\kappa^{top}$
when $Y$ is a sphere.

Now let $Y=(\frac{Z_0(\S^q)}{pZ_0(\S^q)})^0=K(\Z_p; q)$. Let $i:
\S^q \longrightarrow Y$ denote the generator of $\pi_q(Y)$. Let
$u\in H^q(Y; \Z_p)$ and $\sigma\in H_q(Y; \Z_p)$ be given. Since
$i_*:H_q(\S^q; \Z_p)\longrightarrow H_q(Y; \Z_p)$ is an isomorphism,
there is an element $\tau\in H_q(\S^q; \Z_p)$ with $i_*\tau=\sigma$.
Then by equation \ref{naturality of kronecker pairing} and the case
of spheres, $\kappa^{top}(u, \sigma)=\kappa^{top}(i^*u,
\tau)=\kappa(i^*u, \tau)=\kappa(u, \sigma)$.

For a general space $Y$, fix $u\in H^q(Y; \Z_p)$ and let $f:Y
\longrightarrow K(\Z_p, q)$ be the map classifying $u$; i.e.,
$u=f^*\iota$ where $\iota \in H^q(K(\Z_p, q); \Z_p)$ is the
canonical generator. Then for any $\tau\in H_q(Y; \Z_p)$, we have
$$\kappa^{top}(u, \tau)=\kappa^{top}(\iota, f_*\tau)=\kappa(\iota,
f_*\tau)=\kappa(u, \tau).$$ This complete the proofs.
\end{proof}

\section{The Compatibility Of The Duality Theorem With The
$\Z_2$-Poincar\'e Duality} Let us recall that a projective variety
$X$ is said to have full real points if $dim_{\R}ReX=dim_{\C}X$.

\begin{theorem}
Suppose that $X$ is a real projective manifold of dimension $m$ with
full real points and $ReX$ is connected. Then the duality theorem
between reduced real morphic cohomology and reduced real Lawson
homology is compatible with the $\Z_2$-Poincar\'e duality, i.e., the
following diagram commutes:
$$\xymatrix{RL^tH^k(X) \ar[r]^-{\mathscr{RD}} \ar[d]_{\Phi} & RL_{m-t}H_{m-k}(X)
\ar[d]^{\Phi}\\
H^k(ReX; \Z_2) \ar[r]^-{\mathscr{P}} & H_{m-k}(ReX; \Z_2)}$$ where
$\mathscr{P}$ is the Poincar\'e duality map. \label{compatible}
\end{theorem}

\begin{proof}
By Proposition \ref{s map commutes}, we reduce the problem to the
case $t=m$, i.e., we need to prove the commutativity of the
following diagram:
$$\xymatrix{RL^mH^k(X) \ar[r]^{\mathscr{RD}} \ar[d]_{\Phi} & RL_0H_{m-k}(X)
\ar[d]^{\tau}\\
H^k(ReX; \Z_2) \ar[r]^{\mathscr{P}} & H_{m-k}(ReX; \Z_2)}$$ where
$\tau$ is the Dold-Thom isomorphism.

Since the evaluation and the intersection products:
$$H^k(ReX; \Z_2)\otimes H_k(ReX; \Z_2) \overset{<,>}{\longrightarrow} H_0(ReX;
\Z_2),$$
$$H_{m-k}(ReX; \Z_2)\otimes H_k(ReX; \Z_2)
\overset{\odot}{\longrightarrow} H_0(ReX; \Z_2)$$ are perfect
pairings, it suffices to prove that
$$<\Phi(\alpha), \gamma>=\tau(\mathscr{RD}(\alpha))\odot \gamma$$
for all $\alpha\in RL^mH^k(X)$ and all $\gamma\in H_k(ReX; \Z_2)$ .
To prove this equality, it suffices to prove the commutativity of
the following diagram:

\begin{equation}
\xymatrix{H^k(ReX;\Z_2)\otimes H_k(ReX; \Z_2) \ar[rr]^-{<,>} & &
H_0(ReX; \Z_2) \\
\pi_{m-k}(R^m(X))\otimes \pi_k(R_0(X)) \ar[u]^{\Phi\otimes \tau}
\ar[d]_{\mathscr{RD} \otimes pr^*_1}
\ar[rr]^-{\backslash } & &\pi_{m}(R_0(\A^m)) \ar[u]_{\epsilon^{-1}} \ar[d]^{=}\\
\pi_{m-k}(R_m(X\times \A^m))\otimes \pi_k(R_m(X\times \A^m))
\ar[rr]^-{pr_2\circ \bullet} & &\pi_m(R_0(\A^m))\\
H_{m-k}(ReX; \Z_2)\otimes H_k(ReX; \Z_2) \ar[u]^{pr^*_1\circ
\tau^{-1}\otimes pr^*_1\circ \tau^{-1}}  \ar[rr]^-{\odot} & &
H_0(ReX; \Z_2) \ar[u]_{\epsilon}} \label{duality diagram}
\end{equation}

where
$$\epsilon=pr_{2*}\circ pr^*_1: \pi_0(R_0(X))\longrightarrow \pi_m(R_0(X\times \A^m))
\longrightarrow \pi_m(R_0(\A^m))$$ and $\backslash$ is the slant
product. The commutativity of the top square follows from the
construction of slant product \ref{slant product} and the naturality
of Kronecker pairing \ref{kronecker pairing}.

To verify the commutativity of the middle square of the diagram
\ref{duality diagram}, it suffices to show that the following
diagram commutes:
$$
\xymatrix{R^m(X)\wedge R_0(X) \ar[rr]^-{\backslash}
\ar[d]_{\mathscr{RD}\wedge pr^*_1}
& & R_0(\A^m) \ar[d]^{=}\\
R_m(X\times \A^m)\wedge R_m(X\times \A^m) \ar[rr]^-{pr_{2*}\circ
\bullet} & & R_0(\A^m)}
$$

For $f+R_0(\P^{m-1})(X)\in R^m(X)$, $\sum n_ix_i +Z_0(X)^{av} \in
R_0(X)$, $\mathscr{RD}(f+R_0(\P^{m-1})(X))=f+R_m(X\times \P^{m-1})$,
$pr^*_1(\sum n_ix_i+Z_0(X)^{av})=\sum n_i(x_i\times
\P^m)+R_m(X\times \P^{m-1})$, $pr_{2*}[(f+R_m(X\times
\P^{m-1}))\bullet (\sum n_i(x_i\times \P^m)+R_m(X\times
\P^{m-1}))]=pr_{2*}(f\bullet \sum n_i(x_i\times
\P^m))+R_0(\P^{m-1})$ which is the slant product of $f+R_m(X\times
\P^{m-1})$ and $\sum n_ix_i + Z_0(X)^{av}$.

To prove the commutativity of the bottom square, it suffices to
prove the commutativity of the following diagram:
\begin{equation}
\xymatrix{H_{m-k}(ReX; \Z_2)\otimes H_k(ReX; \Z_2) \ar[d]_{\times}
\ar[rr]^-{(pr^*_1)^{\otimes 2}} & & H^{BM}_{2m-k}(ReX\times \R^m;
\Z_2)\otimes H^{BM}_{m+k}(ReX\times \R^m; \Z_2)
\ar[d]^{\times}\\
H_m(ReX\times ReX; \Z_2) \ar[d]_{\Delta!} \ar[rr]^{(pr_1\times
pr_1)^*}
& & H^{BM}_{3m}((ReX\times \R^m)^2; \Z_2) \ar[d]^{\Delta!}\\
H_0(ReX; \Z_2) \ar[rr]^{pr^*_1} \ar[d]_{=} & & H^{BM}_{m}(Re X\times
\R^m; \Z_2)
\ar[d]^{pr_{2*}}\\
H_0(ReX; \Z_2) \ar[rr]^{\epsilon} & & H^{BM}_{m}(\R^m; \Z_2)}
\label{last square}
\end{equation}

The composition of the maps in the right column can be identified
with the map $pr_{2*}\circ \bullet$ in diagram ~\ref{duality
diagram} using the naturality of $\tau$ and the homotopy property of
trivial bundle projection of reduced real Lawson homology and the
composition of the maps in the left column of diagram ~\ref{last
square} is the intersection pairing. Thus the commutativity of
diagram \ref{last square} implies the commutativity of the last
square of \ref{duality diagram}.

The evident intertwining of the external product $\times$ and the
flat pull-back $pr^*_1$ implies the commutativity of the top square.
The Gysin maps and flat pull-backs commute, for a proof, for example
in \cite[3.4.d]{FG}. These maps are real hence induce maps in
reduced real Lawson homology and commute with the Dold-Thom
isomorphism. This give the commutativity of the middle square. The
commutativity of the bottom square comes from the definition of
$\epsilon$.
\end{proof}

From the commutative diagram in Theorem \ref{compatible}, we have
the following result.
\begin{corollary}
If $X$ is a real projective variety of dimension $m$ with full real
points and $ReX$ is connected, then
$$RL^tH^k(X)\cong H^k(ReX; \Z_2), \mbox{ for } 0\leq t \leq m$$
\end{corollary}

\appendix
\section{}
We prove that each cycle group that we deal with in this paper has
the homotopy type of a CW-complex. The proof is inspired by a result
of Friedlander and Walker \cite[Proposition 2.5]{FW}. Throughout
this section, $X, Y$ are projective or real projective varieties of
dimension $m$ and $n$ respectively. To simplify our notation, we
write $\mathscr{C}_d=\mathscr{M}or(X, \mathscr{C}_{r, d}(X))$, the
set of morphisms from $X$ to the Chow variety $\mathscr{C}_{r,
d}(Y)$, and $\mathscr{C}=\coprod_{d\geq 0}\mathscr{C}_{d}$ the
monoid of $Y$-valued $r$-cocycles on $X$. By \cite[Lemma 1.4]{FL1},
$\mathscr{C}_d$ is a quasi-projective variety. Let
$G_n=\coprod_{d_1+d_2\leq n}\mathscr{C}_{d_1}\times
\mathscr{C}_{d_2}$, and $K_n=G_n/\sim$ where $(\alpha_1,
\beta_1)\sim (\alpha_2, \beta_2)$ if and only if
$\alpha_1+\beta_2=\alpha_2+\beta_1$. The map
$$\widetilde{i}_n: K_n \longrightarrow K_{n+1}$$ induced by the inclusion
$i_n: G_n \rightarrow G_{n+1}$ is an inclusion for each $n$. Let
$p_n:G_n \rightarrow K_n$ be the quotient for $n\in \N \cup\{0\}$.
Let $G'_n=p^{-1}_{n+1}(\widetilde{i}_n(k_n))$,
$S_n=\coprod_{a+b=n}\mathscr{C}_a\times \mathscr{C}_b$, and
$R_n=Image\{\coprod_{d_1+d_2+2e=n}\mathscr{C}_{d_1}\times
\mathscr{C}_{d_2}\times \mathscr{C}_{e} \rightarrow S_n\}$ where the
map is given by $(C_1, C_2, C_3) \mapsto (C_1+C_3, C_2+C_3)$. By the
proof of \cite[Proposition 1.3]{FG}, $R_n$ is a subcomplex of $S_n$.
Since $G'_n=G_n\cup R_n$, hence $G'_n$ is a CW-complex. Let
$\widetilde{K}_n=\widetilde{i}_n(K_n)$.

\begin{lemma}
The diagram
$$\xymatrix{G'_n \ar[r]^{i'_n} \ar[d]_{p_{n+1}} & G_{n+1} \ar[d]^{p_{n+1}}\\
\widetilde{K}_n \ar[r]^{\widetilde{i'}_n} & K_{n+1} }$$ is a diagram
of pushout.
\end{lemma}

\begin{proof}
Suppose that the following is a commutative diagram:
$$\xymatrix{G'_n \ar[r]^{i'_n} \ar[d]_{p_{n+1}} & G_{n+1} \ar[d]^f\\
\widetilde{K}_n \ar[r]^g & A }$$  We define a map $\varphi: K_{n+1}
\rightarrow A$ by $\varphi(x)=f(p^{-1}_{n+1}(x))$. Since for $y\in
p^{-1}_{n+1}(x)$ where $x\in \widetilde{K}_n$,
$f(y)=g(p_{n+1}(y))=g(x)$, and if $x\in K_{n+1}-\widetilde{K}_n$,
$|p^{-1}_{n+1}(x)|=1$, hence $f$ is constant on each fibre of
$p_{n+1}$. Therefore $\varphi$ is continuous and $f=\varphi \circ
p_{n+1}$. It is easy to see that $\varphi \circ \widetilde{i'}_n=g$.

Suppose that $\varphi': K_{n+1} \rightarrow A$ is another map which
makes the following diagrams commute:
$$\xymatrix{G'_n \ar[r]^{i'_n} \ar[d]_{p_{n+1}} & G_{n+1} \ar[d]_{p_{n+1}} \ar[rrdd]^{f} &&\\
\widetilde{K}_n \ar[r]^{\widetilde{i'}_n} \ar[rrrd]_{g} & K_{n+1} \ar[rrd]^{\varphi, \varphi'} &&\\
& & & A}$$

If $x\in \widetilde{K}_n$, $\varphi'(x)=g(x)=\varphi(x)$; if $x\in
K_{n+1}-\widetilde{K}_n$, there is a unique $y\in G_{n+1}-G'_n$ such
that $\varphi'(x)=f(y)=\varphi(x)$. Therefore $\varphi=\varphi'$.
\end{proof}

\begin{theorem}
The cycle space $Z_r(Y)(X)$ is a CW-complex. In particular, $Z_r(X)$
is a CW-complex.
\end{theorem}

\begin{proof}
We prove by induction. $K_0=\{0\}$ is a point. Assume that $K_n$ is
a CW-complex and the quotient map from $G_n$ to $K_n$ is regular
cellular. Since $\widetilde{K}_n$ is homeomorphic to $K_n$, and
$G_n, G'_n$ are CW-complexes, hence there is a CW-complex $W_{n+1}$
such that the following diagram is a pushout:
$$\xymatrix{ G'_n \ar[r] \ar[d] & G_{n+1} \ar[d]\\
K'_n \ar[r] & W_{n+1}}$$ and the map from $K'_n$ to $W_{n+1}$ is
regular cellular. Then $W_{n+1}$ is homeomorphic to $K_{n+1}$, hence
$K_{n+1}$ is a CW-complex. And the inclusion map from $K_n$ to
$K_{n+1}$, the quotient map from $G_n$ to $K_n$ are regular
cellular. Therefore $Z_r(Y)(X)$ is a CW-complex.
\end{proof}

As proved by Hironaka in \cite{H} that any projective variety admits
a triangulation by semi-algebraic simplices which can be chosen so
that any specified finite collection of semi-algebraic closed
subsets consists of sub-complexes. Recall that a subset in $\R^n$
defined by the zero loci of some real polynomials is called a
totally real algebraic variety (see \cite{Teh2}). We may embed the
projective space $\P^n$ into $\R^N$ as a totally real algebraic
variety (\cite[Proposition 3.4.6]{BCR}), and the conjugation of
$\P^n$ can be realized as a totally real algebraic map.

\begin{corollary}
The $Y$-valued real cycle space $Z_r(Y)(X)_{\R}$ is a sub-complex of
$Z_r(Y)(X)$. In particular, the real cycle group $Z_r(X)_{\R}$ is a
sub-complex of $Z_r(X)$.
\end{corollary}

\begin{proof}
The set $\mathscr{C}_{d}(Y)(X)_{\R}$ is the set of real points of
$\mathscr{C}_{d}$ which is a totally real quasiprojective variety
and thus $\mathscr{C}_{d}(Y)(X)_{\R}$ is a sub-complex of
$\mathscr{C}_{d}$. Let
$G_{n\R}:=\coprod_{d_1+d_2=n}\mathscr{C}_{d_1}(Y)(X)_{\R}\times
\mathscr{C}_{d_2}(Y)(X)_{\R}$. Since the quotient map $p_n: G_n
\rightarrow K_n$ is regular cellular, hence
$K_n(X)_{\R}:=p_n(G_{n\R})$ is a sub-complex of $K_n$. Therefore
$Z_r(Y)(X)_{\R}$ is a sub-complex of $Z_r(Y)(X)$.
\end{proof}

\begin{corollary}
The $Y$-valued averaged cycle space $Z_r(Y)(X)^{av}$ is a
sub-complex of $Z_r(Y)(X)_{\R}$. In particular, the averaged cycle
group $Z_r(X)^{av}$ is a sub-complex of $Z_r(X)_{\R}$.
\end{corollary}

\begin{proof}
Let $\phi: \mathscr{C}_{d}(Y)(X) \rightarrow \mathscr{C}_{2d}(Y)(X)$
be the map defined by $\phi(c):=c+\overline{c}$. Since the cycle
addition is an algebraic map, and the conjugation is semi-algebraic,
hence $\phi$ is semi-algebraic. Let $\mathscr{C}_{2d}(Y)(X)^{av}$ be
the image of $\phi$, then $\mathscr{C}_{2d}(Y)(X)^{av}$ is a
semi-algebraic set. We may take it as a sub-complex of
$\mathscr{C}_{2d}(Y)(X)_{\R}$. Let
$G^{av}_{2n}=\coprod_{d_1+d_2=n}\mathscr{C}_{2d_1}(Y)(X)^{av}\times
\mathscr{C}_{2d_2}(X)^{av}$. Then $G^{av}_{2n}$ is a sub-complex of
$G_{2n\R}$. Therefore $Z_r(X)^{av}$ is a sub-complex of
$Z_r(X)_{\R}$.
\end{proof}

\begin{theorem}
The $Y$-valued reduced real cycle group $R_r(Y)(X)$ is a CW-complex.
In particular, the reduced real cycle group $R_r(X)$ is a
CW-complex.
\end{theorem}

\begin{proof}
Let $M_n=K_n(Y)(X)_{\R}+Z_r(Y)(X)^{av}$. Since $K_n(Y)(X)_{\R}$ and
$Z_r(Y)(X)^{av}$ are sub-complexes of $Z_r(Y)(X)_{\R}$, and the
cycle addition $+$ is regular cellular, hence $M_n$ is a sub-complex
of $Z_r(Y)(X)_{\R}$. We have an obvious pushout diagram
$$\xymatrix{M_n \ar[r]^{i} \ar[d]_{p_n} & M_{n+1} \ar[d]^{p_{n+1}} \\
K_n(Y)(X)_{\R}+Z_r(Y)(X)^{av} \ar[r] \ar[r] &
K_{n+1}(Y)(X)_{\R}+Z_r(Y)(X)^{av}}$$ By induction,
$K_n(Y)(X)_{\R}+Z_r(Y)(X)^{av}$ is a CW-complex, hence
$K_{n+1}(Y)(X)_{\R}+Z_r(Y)(X)^{av}$ is homeomorphic to the
adjunction space $M_{n+1}\coprod_{p_{n+1}}
K_n(Y)(X)_{\R}+Z_r(Y)(X)^{av}$ consequently, a CW-complex, and the
map from $K_n(Y)(X)_{\R}+Z_r(Y)(X)^{av}$ to
$K_{n+1}(Y)(X)_{\R}+Z_r(Y)(X)^{av}$ is regular cellular, therefore,
$R_r(Y)(X)$ is a CW-complex.
\end{proof}

\begin{theorem}
The reduced real cocycle $R^t(X)$ is a CW-complex.
\end{theorem}

\begin{proof}
Let $K_n=\coprod_{d_1+d_2\leq n}\mathscr{C}_{d_1}(\P^t)(X)\times
\mathscr{C}_{d_2}(\P^t)(X)/\sim$,
$M_n=K_n+Z_0(\P^t)(X)^{av}+R_0(\P^{t-1})$, $N_n=p_n(M_n)$ where
$p_n$ is the restriction of the quotient map $p: R_0(\P^t)(X)
\rightarrow R^t(X)$ to $M_n$. Then we have a pushout diagram
$$\xymatrix{M_n \ar[r] \ar[d]_{p_n} & M_{n+1} \ar[d]^{p_{n+1}}\\
N_n \ar[r] & N_{n+1}}$$ By induction, we show that $R^t(X)$ is a
CW-complex.
\end{proof}

\begin{acknowledgement}
The author thanks Blaine Lawson for his encouragement and guidance
of his study, and for his remarks and corrections of this paper. He
also thanks Christian Haesemeyer for his thorough proofreading and
the National Center for Theoretical Sciences in Hsinchu, Taiwan for
its warm hospitality.
\end{acknowledgement}

\begin{affiliation}
   Jyh-Haur Teh\\
   National Tsing Hua University,\\
Department of Mathematics,\\ Third General Buidling,\\
No. 101, Sec 2, Kuang Fu Road, Hsinchu 30043, Taiwan.\\
   Email: jyhhaur@math.nthu.edu.tw\\

\end{affiliation}

\end{document}